\documentstyle[12pt]{report}

\newcommand\A{{\cal A}}
\renewcommand\a{\alpha}

\renewcommand\c{\circ}
\newcommand\C{{\mbox{\rm\bf C\hspace{-7.3pt}}{^{_{\bf\mid}}}\hspace{4.5pt}}}
\newcommand\CC{{\cal C}}
\renewcommand\d{\dots}
\newcommand\D{\Delta}
\newcommand\e{\eta}
\newcommand\E{{\cal E}}
\newcommand\g{\gamma}
\newcommand\G{{\cal G}}
\newcommand\Gl{\op{Gl}}

\newcommand\J{{\cal J}}
\newcommand\k{\kappa}

\renewcommand\ll{\lambda}
\newcommand\m{\mu}
\newcommand\n{{\bf\Im}{\hspace{-1.5pt}{\bf\wr}\hspace{-1pt}}}
\newcommand\NN{N^{(2)}}
\newcommand\N{{\cal N}}

\renewcommand\O{\Omega}
\newcommand\oo{\omega}
\newcommand\ok{{\cal O}}
\newcommand\ol{\ok_L(x)}
\newcommand\om{\ok_M(y)}
\newcommand\op[1]{\mathop{\rm #1}\nolimits}
\newcommand\p{\partial}
\newcommand\po{$\!\!\!{\bf .}$ }
\renewcommand\P{\Phi}
\newcommand\PP{\Phi^{(2)}}
\newcommand\r{PH}
\newcommand\R{{\rm I\hspace{-2.5pt} R}}
\renewcommand\t{\times}
\newcommand\te{\theta}

\newcommand\T{\Theta}
\newcommand\U{\Upsilon}
\renewcommand\u{{\cal U}}
\newcommand\ve{\varepsilon}
\newcommand\vp{\varphi}
\newcommand\V{{\cal V}}
\newcommand\w{{\cal W}}
\newcommand\x{\xi}
\newcommand\Y{{\cal Y}}
\newcommand\z{\sigma}
\newcommand\Z{{\rm Z\mkern-5muZ}}

\def\Rom#1{\uppercase\expandafter{\romannumeral#1}}

\newcommand\1{{\bf 1}}
\newcommand\2{{\wedge2}}
\newcommand\3{{\vec\te}^{k-2}}

\newcommand\qed{\phantom{\underline{y}}\hfill\hfill$\Box$}

\newcommand\lv{\vert\![}
\newcommand\rv{]\!\vert}
\newcommand\phd{\phantom{\dfrac12}}
\newcommand\bib[1]{\bibitem[#1]{#1}}

\newcommand{\text}[1]{{\mbox{\rm #1}}}
\newcommand{\dfrac}[2]{\frac{\displaystyle #1}{\displaystyle #2}}
\newcommand\add{\addtocounter{a}{1}}

\newcommand{\olra}[1]{\stackrel{#1}{\longrightarrow}}
\newcommand{\barwedge}{\bar{\wedge}}

\newtheorem{th}{Theorem}
\newtheorem{prop}{Proposition}
\newtheorem{lem}{Lemma}

\newenvironment{ex}[1]{\trivlist \item[\hskip \labelsep{\bf Example #1.}]}%
{\endtrivlist}
\newenvironment{cor}{\trivlist \item[\hskip \labelsep{\bf Corollary.}]}%
{\endtrivlist}
\newenvironment{dfn}[1]{\trivlist \item[\hskip \labelsep{\bf Definition #1.}]}%
{\endtrivlist}
\newenvironment{lll}[1]{\trivlist \item[\hskip
\labelsep{{\bf Lemma} ${\bf #1}_{\bf k}$.}]}%
{\endtrivlist}
\newenvironment{rk}[1]{\trivlist \item[\hskip
\labelsep{{\it\underline{Remark #1}.\/}}]}%
{\endtrivlist}
\newenvironment{proof}{\trivlist \item[\hskip
\labelsep{{\it\underline{Proof}.\/}}]}%
{\endtrivlist}
\newenvironment{Proof}[1]{\trivlist \item[\hskip
\labelsep{{\it\underline{#1}.\/}}]}%
{\endtrivlist}

\newcounter{a}
\newcounter{f}
\makeatletter
\newcommand{\@thefnmark}{$^\fnsymbol{f}$}
\renewcommand{\@makefnmark}{\hbox{\mathsurround=0pt
                           $^{\fnsymbol{f}}$}}
\renewcommand{\@makefntext}[1]{\parindent=1em\noindent
            \hbox to 1.8em{\hss$^{\fnsymbol{f}}$}#1}
\makeatother

\begin{document}

\title{\bf Nijenhuis tensors and obstructions for pseudoholomorphic mapping
constructions}
%\rightmark{Nijenhuis tensors and PH-mapping}

\author{\bf Boris~S.~Kruglikov}

\date{October 1, 1996}

\begin{abstract}

In this paper we present some approaches to classification of almost
complex structures and to construction of local or formal
pseudoholomorphic mapping from one almost complex manifold to another.
The corresponding criteria are given in terms of Nijenhuis tensors and
their generalizations.

We deal with the prolongations of the Cauchy-Riemann equation in 1-jets of
the mapping from one almost complex manifold to another. We give a criterion
of the prolongation existence of $k$-th pseudoholomorphic jets to the
$(k+1)$-th ones. As a consequence
basing on the Kuranishi's approach
we obtain formal normal forms of almost complex structures. As another
consequence we get the formal part of Newlander-Nirenberg's theorem on
integrability.

We also introduce the notion and study the properties of the linear
Nijenhuis tensor, which we apply to almost complex structures of general
position. Two applications are exhibited: a classification of
four-dimensional almost complex structures in terms of distributions and
a characterization of the new class of structures which we call Lie
almost complex structures.

\end{abstract}

\maketitle
%\thanks

\tableofcontents

\clearpage

%%%%%%%%%%%%%%%%%%%%%%%%%%%%%%%%%%%%%%%%%%%%%%%%%%%%%%%%%%%%%%%%%%%%%%%%%%

\chapter*{Introduction}
\addcontentsline{toc}{chapter}{\bf\quad \  Introduction}

\hspace{13.5pt}
Let $(L^{2l},j_L)$ and $(M^{2m},j_M)$ be almost complex manifolds,
i.e. manifolds equip\-ped with automorphisms $j_L$ and $j_M$ of the tangent
bundles which satisfy the conditions
$j_L^2=-\1_L: T_*L\to T_*L$ and $j_M^2=-\1_M$. A pseudoholomorphic
(PH-)mapping $u:L\to M$ is such a mapping that its differential
$\left.u_*\right\vert_x=du_x:T_xL\to T_{u(x)}M$ commutes with the almost
complex structure: $j_M\c u_*=u_*\c j_L$.

If $M=\C$ a (nontrivial) pseudoholomorphic mapping $u: L\to M$ is
called a complex
coordinate and in general case there is no such a coordinate
(moreover if on a manifold $L^{2l}$ there exist $l$ independent local
complex coordinates then it is complex). For $l=1$, i.e. when
$(L,j_L)$ is a Riemann surface, PH-mappings do exist and moreover the
family of such mappings possesses the manifold structure (under some
additional conditions, see~\cite{Gr1}, \cite{MS}). If $l>1$ then
for an almost complex structure $j_L$ of general position
such a mapping does not exist even locally. An obstruction is the Nijenhuis
tensor. It is well-known (\cite{NN}, \cite{NW}) that if the Nijenhuis tensor
equals zero identically then the almost complex manifold is complex and for
a pair of such manifolds $L$ and $M$ there exists a local holomorphic
mapping $u:(L,x)\to(M,y)$ (for more see~\cite{K}).
In this paper we consider the existence conditions
of formal and local pseudoholomorphic mappings for general almost complex
structures.

It is worth mentioning that the existence of a (nontrivial)
pseudoholomorphic mapping is a
rare situation for mappings $u:L\to M$ between almost complex manifolds
$(L^{2l},j_L)$ and $(M^{2m},j_M)$ in general position unless $l=1$. For
example, for any $N\ge1$ the pseudoholomorphic imbedding to the standard
complex space $(L^{2l},j_L)\to(\C^N,j_0)$ exists if and only if the structure
$j_L$ is integrable, i.e. complex, which is of course not the case of general
position. This rigidity result exhibits a great difference between almost
complex structures and some other tensor fields such as Riemannian metrics.
Recall that according to the Nash's theorem~\cite{N} any Riemannian manifold
admits an isometric imbedding $(V,g)\to(\R^N,g_0)$ to the Euclidean space
of some big dimension $N\gg1$ (see also~\cite{Gr2}). Thus any
pseudoholomorphic mapping $(L^{2l},j_L)\to(M^{2m},j_M)$, $l>1$, is very
unstable under the perturbation of almost complex structures, and hence due
to the strict rigidity of local pseudoholomorphic mappings the problem
generally goes onto the formal level.

The structure of the paper is the following.
In chapter~1 we prove a necessary and sufficient condition for the
existence of a prolongation of a pseudoholomorphic on 1-jet level mapping
of one almost complex manifold to another to a pseudoholomorphic mapping
on 2-jet level.
In other words we describe the first prolongation of the corresponding
Cauchy-Riemann equation on the 1st-jet level. In this case the
Cauchy-Riemann equation is the commutation condition of the required
mapping $u$ differential $\Phi=du$ at a point with the almost complex
structures: $j_M\c\Phi=\Phi\c j_L$. We introduce the space of
$k$-jets of pseudoholomorphic mappings which does not have the manifold
structure in general case but does have for formally integrable "completed"
Cauchy-Riemann equations. Here by "completion" we mean an addition of
defining relations, i.e. a contraction from the point of view of jet spaces
(\cite{KLV}). The Cauchy-Riemann equation for mappings of one almost
complex manifold to another is formally integrable in the integrable case
only, i.e. when the manifolds are complex. In general case this equation
must be "completed". The first term of "completion" is the commutation
condition of the required mapping differential with the Nijenhuis tensors:
$N_{j_M}\c\wedge^2\Phi=\Phi\c N_{j_L}$. In chapter~2 we consider the
question of formal integrability, introduce a new invariant --- Nijenhuis
tensor of the second order $\NN_j$ and write down the next term of
"completion": $\NN_{j_M}\c\wedge^2(\wedge^2\P)=\P\c\NN_{j_L}$.

In chapter 3 we prove the solvability criterion for the Cauchy-Riemann
equation on the higher jet spaces. As a corollary we get the formal part of
the well-known Newlander-Nirenberg theorem on the integrability of almost
complex structures. We also use this criterion to find formal normal forms of
almost complex structures.

In chapter 4 we introduce the space of linear Nijenhuis tensors and
consider the recovery question for liner complex structure by a Nijenhuis
tensor. In chapter 5 we define the notion of the Nijenhuis tensors of general
position. We call an almost complex structure the structure of general
positions if the corresponding Nijenhuis tensor is in general position.
For the manifold equip\-ped with such structures the necessary and sufficient
condition for existence of a pseudoholomorphic mapping may be formulated in
terms of Nijenhuis tensors only, which forms the statement of theorem~6.

In section 6.1 we present a classification of almost complex structures of
general positions on four-dimensional manifolds in terms of
distributions. In section 6.2 we introduce the notion of Lie almost complex
structure and present a classification of such structures of general
position.

In appendix we prove that the notion general position for the Nijenhuis
tensors, introduced in chapter~3, actually satisfies the general position
properties.

Chapters 4--6 do not use the results and methods of chapters 1--3 and might
be read independently.

The author is grateful to prof. V.\,V.~Lychagin for a warm attention to the
work and helpful discussions.

%%%%%%%%%%%%%%%%%%%%%%%%%%%%%%%%%%%%%%%%%%%%%%%%%%%%%%%%%%%%%%%%%%%%%%%%%%
% 1 %
\chapter{Nijenhuis tensor as the first obstruction \protect\\
to the construction of pseudoholomorphic mapping}
\markboth{1. Nijenhuis tensor as the first obstruction}
{1. Nijenhuis tensor as the first obstruction}

\hspace{13.5pt}
For the mapping $u:(L,x)\to(M,y)$ we denote by $u_*^{\wedge k}: \wedge^k
T_xL\to\wedge^k T_yM$ the induced mapping.

Let us denote by $\m=\m_x$ the ideal in the ring of functions
$C^\infty(L)$ vanishing at the point $x$. The power of this ideal $\m^k$
is the ideal in the ring of functions vanishing
at the point $x$ together with their derivatives of the order $<k$:
$\m^k=\{f \vert \p^\z f(x)=0,\ |\z|<k\}$. Let's denote by $\m^k=\m^k T(p,q)$
the submodule in the module of tensors of the type $(p,q)$ vanishing at the
point $x$ up to the order $k$. Denote by the symbol $\ok_L(x)$
the germ of neighborhoods of the point $x\in L$.
Let us also denote $\mu^k_{x,y}=\{\P: T\ok_L(x)\to T\ok_M(y) |
\op{Im}\P\subset\m^k_y\}$.

A formal pseudoholomorphic mapping is a mapping
$u: \ok_L(x)\to\ok_M(y)$ such that $j_M\c u_*-u_*\c j_L\in\mu^\infty$.
The mapping is called nontrivial if $u_*\notin\mu^\infty$.
Moreover we will consider, in general, nondegenerate mappings:
$u_*\notin\m$ or $du_x\ne0$.

Let us consider the Nijenhuis tensors $N_{j_L}$ and $N_{j_M}$ on manifolds
$L$ and $M$:
 $$
N_j(\x,\e)=[j\x,j\e]-j[j\x,\e]-j[\x,j\e]-[\x,\e].
 $$
Immediately from the definition it follows that if
$u: L\to M$ is a \r-mapping then $N_{j_M}\c u_*^{\wedge2}=u_*\c N_{j_L}$.
Moreover this equality holds at the point $x$ if $u$ preserves the almost
complex structure up to the second order:
$j_M\c u_*-u_*\c j_L\in\mu^2_{x,y}$.

Let us denote $J^0_{PH}(L,M)=J^0(L,M)=L\t M$ (we will use the notations
from the jet spaces theory, see \cite{KLV}, \cite{KS}).
Let us denote by $J^1_{PH}(L,M)$ the subbundle of the jet bundle
$J^1(L,M)\to J^0(L,M)$ with the fiber
 $$
(J^1_{PH})_{x,y}=\{\P: T_xL\to T_yM \,\vert\, j_M\c\P=\P\c j_L\}.
 $$
In other words the manifolds $J^1_{PH}$ is a bundle over $L\t M$
with the fiber $(J^1_{PH})_{x,y}$ being the factor space of the space
$\{u_*: T\ol\to T\om \,| j_M\c{u_*}-u_*\c j_L\in\m_{x,y}\}$ by the space
$\{u_*\in\m_{x,y}\}$. Let us also denote
 $$
J^k_{PH}\supset (J^k_{PH})_{x,y}= \frac
{\{u_*: T\ol\to T\om \,|\, u_*\c j_L - j_M\c u_* \in\m^k_{x,y}\}}
{\{u_*: T\ol\to T\om \,|\, u_*\in\m^k_{x,y}\}}.
 $$
Unlike $J^1_{PH}$ the set $J^k_{PH}$, $k\ge2$, does not have in general
the manifold structure. This is connected with the fact that for the mapping
$\pi^{r,s}_{PH}$ of the restriction of the projection
$\pi^{r,s}: J^r(L,M)\to J^s(L,M)$, $r\ge s$,
the inverse image might be empty: $(\pi^{r,s}_{PH})^{-1}(pt)=\emptyset$.

We apply the Cartan's method of prolongations-projections (\cite{ALV}) to
obtain a criterion of existence of a formal \r-mapping. Denote by $\E^k$ the
image of the projection
$\op{Im}\pi^{k,1}_{PH}\subset J^1_{PH}\subset J^1(L,M)$.
One may consider these subsets in the 1-jet spaces as differential equations
of the first order. The problem of formal pseudoholomorphic mapping
construction is reduced to the contraction of the equation --- to the finding
of a projective limit
 $$
\E^\infty\subset\d\E^2\subset\E^1\subset J^1(L,M).
 $$
The equation $\E^1=J^1_{PH}$ is fibered over $J^0=L\times M$ with the fiber
 $$
\E^1_{{x,y}}=(\pi^{{1,0}}_{PH})^{-1}(x,y)=
\{\P:T_xL\to T_yM \,|\, j_M\c\P=\P\c j_L\}.
 $$

 \begin{th}\po
The restriction of the projection $\pi^{1,0}_{PH}$ to the equation
$\E^2\subset J^1_{PH}$ has the inverse images
 $$
\E^2_{x,y}=
\{\P: T_xL\to T_yM \,|\, j_M\c\P=\P\c j_L,\ N_{j_M}\c \P^\2=\P\c N_{j_L} \}.
 $$
 \end{th}

{\it Notations.\/} For the proof of this and further theorems we need the
calculations in jet spaces. For this purpose we introduce the following
constructions.

 \begin{dfn}{1}
Let us define for arbitrary vector field
$\x$ on $X$ the mapping
$d^p\x=d^p_\nabla\x: \otimes_1^pT_xX\to T_xX$, $x\in X$ by the formula
 $$
d\x(\e)=\nabla_\e\x,\
d^2\x(\e_1,\e_2)=\nabla_{\e_2}\nabla_{\e_1}(\x)-\nabla_{\nabla_{\e_2}\e_1}
(\x),\ {\text{ and for arbitrary }}p:
 $$
 $$
 \begin{array}{l}
d^{p+1}\x(\e_1,\d,\e_{p+1})\\
\qquad {}=d\left(d^p\x(\e_1,\d,\e_p)\right)(\e_{p+1})-
\sum_{i=1}^p d^p\x(\e_1,\d,d\e_i(\e_{p+1}),\d,\e_p),
 \end{array}
 $$
where $\nabla$ is some symmetric connection. If the curvature of the
connection $\nabla$ vanishes then the differential $d^p_\nabla\x$ is
symmetric, $d^p_\nabla\x:
S^pT_xX\to T_xX$, and in local Euclidean coordinates it has the form:
 $$
d^p\x(\e_1,\d,\e_p)^i=\frac{\p^p\x^i}{\p x^{i_1}\d \p x^{i_p}}
\e_1^{i_1}\d\e_p^{i_p},
 $$
and it's easy to see that for $d\x=d^1\x$ holds
$L_\x\e=[\x,\e]=d\e(\x)-d\x(\e)$.
 \end{dfn}

 \begin{dfn}{2}
For arbitrary morphism of the bundles $A=(A_x)_{x\in X}$,
$A_x: \wedge^qT_xX\to T_xX$, let us define the mapping
$d^pA$, $d^pA_x: \wedge^qT_xX\otimes(\otimes_1^pT_xX)\to T_xX$
by induction:
\quad $d^0A:=A$, \quad $d^{p+1}A(\x_1,\d,\x_{q+p+1})=$
 $$
=d\left(d^pA(\x_1,\d,\x_{q+p})\right)(\x_{q+p+1})-
\sum_{i=1}^{q+p} d^pA(\x_1,\d,d\x_i(\x_{q+p+1}),\d,\x_{q+p}).
 $$
 \end{dfn}

If the curvature of the connection $\nabla$ vanishes (which always might be
supposed locally true) then the differential
$d^pA_x: \wedge^qT_xX\otimes S^pT_xX\to T_xX$
has the following form in local coordinates:
 $$
\left[d^pA(\x_1,\d,\x_{q+p})\right]^i=
\sum_{i_1<\d<i_q}\frac{\p^pA^i_{i_1\d i_q}}{\p x^{i_{q+1}}\d\p x^{i_{q+p}}}
\x_1^{i_1}\d \x_q^{i_q}\d \x_{q+p}^{i_{q+p}}.
 $$

Thus a Euclidean connection $\nabla$ defines the mapping
 $$
d^p_\nabla: \O^q\otimes {\cal D} \to \O^q\otimes S^p\O^1\otimes {\cal D},
 $$
where ${\cal D}$ is the module of vector fields, and $\O^q$ is the module of
exterior forms. Let us note that such a mapping might be constructed by any
connection $\nabla$ (\cite{P}), but as far as all the arguments below using
$\nabla$ are local we will suppose in general that this connection is trivial
(Euclidean). We will write $d^p$ instead of $d_\nabla^p$.

 \begin{Proof}{Proof of theorem 1}
1-\r-jet $\P$ is lifted to 2-\r-jet iff the following equations hold:
 \begin{equation}
\left\{
  \begin{array}{l}
j_M\c\P(\x)=\P\c j_L(\x), \ \ \ \ \ \ \ \ \ \ \ \ \ \ \ \ \ \ \,\,
\P^{(2)}(\x,\e)=\P^{(2)}(\e,\x),\phantom{\dfrac12}\\
dj_M(\P\x,\P\e)+j_M\c\P^{(2)}(\x,\e)= \P^{(2)}(j_L\x,\e)+\P\c dj_L(\x,\e),
\phantom{\dfrac12}
% \tag 1
\add
  \end{array}
\right.
 \end{equation}
where $\x,\e\in T_xL$ and $\P^{(2)}: S^2T_xL\to T_yM$ is the 2-symbol of
the required \r-mapping (in local coordinates $\P^{(2)}(\x,\e)^i=
\dfrac{\p^2u^i}{\p x^r\p x^s}\x^r\e^s$). Let us consider the tensor
 \begin{equation}
P\in T^*_xL\otimes T^*_xL\otimes T_yM,\
P(\x,\e)=j_M\c\P^{(2)}(\x,\e)-\P^{(2)}(j_L\x,\e).
% \tag 2
\add
 \end{equation}

From (1) we have:
 \begin{equation}
P(\x,\e)=\P\c dj_L(\x,\e)-dj_M(\P\x,\P\e).
% \tag 3
\add
 \end{equation}

% L1 %
 \begin{lem}\po
The tensor $P$ defined by a symmetric tensor $\P^{(2)}$
by the formula (2) satisfies the equations
 \begin{eqnarray}
(i)\ &&
\quad\quad\quad\quad\quad\ \ \
j_M\c P(\x,\e)=-P(j_L\x,\e),
\quad\quad\quad\quad\quad\ \ \ \
% \tag 4
\add                                                  \\
(ii) &&
\quad\quad\quad\quad
P(\x,\e)-P(\e,\x)=P(j_L\x,j_L\e)-P(j_L\e,j_L\x).
\quad\quad\quad\quad
% \tag 5
\add
 \end{eqnarray}
 \end{lem}

% L2 %
 \begin{lem}\po
Every tensor satisfying the system (4)-(5) may be
represented in the form (2) where $\P^{(2)}$ is a symmetric tensor.
 \end{lem}
 \begin{Proof}{Proof of lemma 2}
Since any bilinear map is decomposed to the sum of
%\linebreak
($j_L$-$j_M$-)linear and antilinear by the first argument maps it follows
from (4) that
 \begin{equation}
P(\x,\e)=j_M\c B(\x,\e)-B(j_L\x,\e)
% \tag 6
\add
 \end{equation}
for some $B$. The lemma states that the map $B:T_xL\otimes
T_xL\to T_yM$ might be chosen symmetric. Let us consider the decomposition
of $B$ into the symmetric and antisymmetric components:
 $$
B=B_0+B_1,\ B_i(\x,\e)=\frac12[B(\x,\e)+(-1)^iB(\e,\x)].
 $$
The substitution of this expression into (5) taking into account (6) gives:
 \begin{equation}
j_M\c B_1(\x,\e)-B_1(j_L\x,\e)-B_1(\x,j_L\e)-j_M\c B_1(j_L\x,j_L\e)=0.
% \tag 7
\add
 \end{equation}

Let us decompose $B_1$ into linear and antilinear by the first and the
second argument components:
 \begin{eqnarray*}
&& B_1=B_1^{0,0}+B_1^{0,1}+B_1^{1,0}+B_1^{1,1},\  B_1^{1,1}\equiv0,\quad \\
&& B_1^{0,0}(\x,\e)=\frac12[B_1(\x,\e)-B_1(j_L\x,j_L\e)],\quad\\
&& B_1^{0,1}(\x,\e)=\frac12[B_1(\x,\e)+j_M\c B_1(\x,j_L\e)],\quad\\
&& B_1^{1,0}(\x,\e)=\frac12[B_1(\x,\e)+j_M\c B_1(j_L\x,\e)].\quad
 \end{eqnarray*}
Using (7) one could easily verify the properties
 $$
B_1^{r,s}(\x,\e)=(-1)^{r+1}j_M\c B_1^{r,s}(j_L\x,\e)=
(-1)^{s+1}j_M\c B_1^{r,s}(\x,j_L\e).
 $$
Note that the tensor
 $$
{\hat B}_0(\x,\e)=B_1^{0,1}(\x,\e)-B_1^{1,0}(\x,\e)=\frac{j_M}2
[B_1(\x,j_L\e)-B_1(j_L\x,\e)]
 $$
is symmetric and the tensor $B_1+{\hat B}_0=B_1^{0,0}+2B_1^{0,1}$ is
($j_L$-$j_M$-)linear by $\x$. Therefore the transformation
 $$
B=B_0+B_1\mapsto B-(B_1+{\hat B}_0)=B_0-{\hat B}_0
 $$
determines the required symmetric tensor $\P^{(2)}=B_0-{\hat B}_0$,
which satisfies (2).\qed
 \end{Proof}

Thus the equation $\E^2$ is singled out as a subequation in $\E^1$ by
relations (3)-(5). At that the equation (4) gives the identity according to
 \begin{equation}
j_M\c\P=\P\c j_L,\ j_L\c dj_L=-dj_L\c(j_L\otimes\1),\
j_M\c dj_M=-dj_M\c(j_M\otimes\1).
% \tag 8
\add
 \end{equation}
Under substitution of equation (3) into (5) we get
 $$
 \begin{array}{r}
dj_M(\P\x,\P\e)-dj_M(\P\e,\P\x)-dj_M(j_M\P\x,j_M\P\e)+dj_M(j_M\P\e,j_M\P\x)=
\phantom{\dfrac12}\\
\P\left(dj_L(\x,\e)-dj_L(\e,\x)-dj_L(j_L\x,j_L\e)+dj_L(j_L\e,j_L\x)\right).
\phantom{\dfrac12}
 \end{array}
 $$
Multiplying this equality by $j_M$ and using (8) and the identity
 \begin{equation}
N_j(\x,\e)=-dj(j\x,\e)-dj(\x,j\e)+dj(j\e,\x)+dj(\e,j\x),
% \tag 9
\add
 \end{equation}
we get the equation $N_{j_M}\c\P^\2=\P\c N_{j_L}$, and the theorem is proved.
\hfill$\Box$
 \end{Proof}

Thus the theorem gives the equations at a point (on $T_xL$) determining the
possibility of pseudoholomorphic mapping construction on 2-jet level.

Let us introduce the linear space of Nijenhuis tensors at the complex linear
space $(\C^n,j)$ (further it will be shown in theorem~5 (chapter~4)
that every tensor
from this space can be realized as the Nijenhuis tensor of some almost
complex structure):
 \begin{equation}
\N^{(n)}_j=\{N\in \wedge^2(\C^{n})^*\otimes\C^{n}\ \vert \ N(j\x,\e)=
N(\x,j\e)=-jN(\x,\e)\}.
% \tag {10}
\add
 \end{equation}
Consider on $\N^{(n)}_j$ the factorization mapping
 \begin{equation}
\te: \N_j^{(n)}\to \N_j^{(n)}/\Gl_\C(n),
% \tag {11}
\add
 \end{equation}
which associate to each tensor its orbit under the full group of complex
linear transformations. The Cartan's approach to formal integrability (and
also one of Spencer, see Goldshmidt's theorem in~\cite{KS}, \cite{ALV})
assumes to find by the prolongations-projections method the equation
$\E^{(\infty)}\subset\E^2\subset J^1_{PH}$,
which prolongs to the 2-jet space and the symbols of which form a bundle.
In particular, one needs to require that
$(\E^2)_{x^\prime,y^\prime}\ne0$ for all $x^\prime\in L$,
$y^\prime\in M$. For example, it is so if one requires the constancy of
the natural mappings $\te^{(N)}_L:\ol\to \N_{j_L}^{(l)}/\Gl_\C(l)$,
$\te^{(N)}_M:\om\to \N_{j_M}^{(m)}/\Gl_\C(m)$, associating to a point $x$
the image of the Nijenhuis tensor $(N_j)_x$ at this point upon the mapping
(11), and the existence of a morphism $\P: T_xL\to T_yM$, which conjugates
$N_{j_L}$ with $N_{j_M}$. Manifolds $(L,j_L)$ with
$\te^{(N)}_L\equiv\op{const}$ do exist as the following example shows. Note
also the complexity of the linear classification
problem for linear Nijenhuis tensors $N\in\N^{(n)}_j$.
This problem contains Lie algebras classificational
problem as a subproblem.

 \begin{ex}{1}
Let's consider a Lie group $G$. The Lie algebra
of the Lie group $G\times G$ is isomorphic to $\G\oplus\G$. Let's introduce
the almost complex structure on $T_e(G\times G)$ by formula
$j(\x,\e)=(-\e,\x)$ and extend it by the left shifts onto the whole group.
The Nijenhuis tensor is left-invariant and has the form:
 $$
N((\x,0),(\e,0))=-N((0,\x),(0,\e))=(-[\x,\e],[\x,\e]),
 $$
 $$
\quad\quad\quad\quad\ \ \ \
N((\x,0),(0,\e))=N((0,\x),(\e,0))=([\x,\e],[\x,\e]).
\quad\quad\quad\quad\ \ \ \
\Box
 $$
 \end{ex}

 \begin{rk}{1}
In the paper~\cite{Gau} it is shown that there exists a canonical almost
complex structure $j^{[1]}$ on the manifold $J^1_{PH}(L,M)$. This structure
is defined as $(j_L\oplus j_M)\oplus j_M$ on $T_\bullet J^1_{PH}=
T_\bullet J^0_{PH} \oplus Smbl_{\E^1}^1$ in the decomposition induced by
minimal connections on manifolds $L$ and $M$ (a minimal connection is by
definition such that its torsion form coincides with the torsion form
of the almost complex structure, i.e.
with the corresponding Nijenhuis tensor, see~\cite{L}). At the points
$p\in\E^2\subset J^1_{PH}$ the almost complex structure $j^{[1]}$ on
$T_pJ^1_{PH}$ may be defined easier. Actually, 1-\r-jet $p$ prolongs
to 2-\r-jet $p^{(2)}$, and $p^{(2)}$ determines a subspace $L_{p^{(2)}}$ in
the tangent space to $J^1_{PH}$, i.e. the tangent space to the graphic
inducing $p^{(2)}$, cf.~\cite{KLV}. Thus, the tangent space $T_pJ^1_{PH}$ is
isomorphic to the direct sum of the space $Smbl^1_{\E^1}$ --- the tangent
space to the linear space-fiber and the space $L_{p^{(2)}}$, and since both
summand are complex, there exists a canonical complex structure on
$T_pJ^1_{PH}$. Note that even though the prolongation $p^{(2)}$ may be
chosen not in unique way the structure $j_p^{[1]}$ does not depend on it.
Suppose that the subsets $(\E^2)_{x^\prime,y^\prime}$ do not depend on
points $x^\prime\in L$, $y^\prime\in M$ (compare with the above arguments)
and form a bundle over $J^0_{PH}$, i.e. $\E^2$ is a submanifold in $\E^1$.
It follows from theorem~1 that this submanifold is not, in general case,
complex, and its difference from complex one "measures" the nonintegrability
of the structures $j_L$ and $j_M$. Actually, due to the antilinearity, the
quadric determined in the complex linear space
$\{\Phi\in T_x^*L\otimes T_yM \,\vert\, j_M\c\Phi=\Phi\c j_L\}$ by the
condition $\{ N_{j_M}\c\Phi^\2= \Phi\c N_{j_L} \}$ is not a complex
submanifold iff $\left.(N_{j_M})\right\vert_{\op{Im}\Phi^\2}\not\equiv0$.
Hence on $J^2_{PH}$ does not exist a complex structure such that the
projection restriction $\pi^{2,1}_{PH}: J^2_{PH}\to J^1_{PH}$ is
pseudoholomorphic.

As long as $J^1_{PH}(L,M)$ possesses the canonical almost complex structure
$j^{[1]}_{L,M}$, we can consider the manifold
$J^{[2]}_{PH}(L,M)= J^1_{PH}(L,J^1_{PH}(L,M))$ with the almost complex
structure $j_{L,M}^{[2]}$, and define by induction the manifolds
$J^{[k+1]}_{PH}(L,M)$= $J^1_{PH}(L,J^{[k]}_{PH}(L,M))$ with almost complex
structures $j_{L,M}^{[k+1]}$. Note that
$J^k_{PH}(L,M)\subset J^{[k]}_{PH}(L,M)$.
We say that 1-jet $\Phi\in(J^1_{PH}(L,M))_{x,y}$ is lifted up to the
$k$-jet if $(J^k_{PH})_{x,y}\supset(\pi^{k,1}_{PH})^{-1}(\Phi)
\ne\emptyset$. According to theorem~1, 1-jet $\Phi$ is lifted up to
a 2-jet iff $N_{j_M}\c\Phi^\2=\P\c N_{j_L}$, and in this case there exists
(not necessarily unique) a lift
$\Phi^{[2]}\in (J^2_{PH})_{x,y}\subset (J^{[2]}_{PH})_{x,y}$,
which could be considered as a 1-jet in the ambient space. This 1-jet can
be lifted up to a 2-jet $\Phi^{[3]}\in (J^{[3]}_{PH})_{x,y}$ iff
$N_{j_{L,M}^{[1]}}\c({\Phi^{[2]}})^\2={\P^{[2]}}\c N_{j_L}$.
In general case, if among lifts $\Phi^{[k]}\in (J^{[k]}_{PH})_{x,y}$ of
the jet $\Phi^{[k-1]}\in (J^{[k-1]}_{PH})_{x,y}$ there exists such that
 \begin{equation}
N_{j_{L,M}^{[k-1]}}\c({\Phi^{[k]}})^\2={\P^{[k]}}\c N_{j_L},
% \tag{12}
\add
 \end{equation}
then the jet $\Phi^{[k]}\in (J^{[k]}_{PH})_{x,y}$ may be lifted up to a jet
$\Phi^{[k+1]}\in (J^{[k+1]}_{PH})_{x,y}$. Thus we have a sequence of
necessary conditions (12), and upon their fulfillment there exist lifts
$\Phi^{[r]}\in (J^{[r]}_{PH})_{x,y}$ of 1-jet $\Phi\in (J^1_{PH})_{x,y}$,
and one can try to find such lifts that $\Phi^{[r]}\in (J^r_{PH})_{x,y}$.
\qed
 \end{rk}

%%%%%%%%%%%%%%%%%%%%%%%%%%%%%%%%%%%%%%%%%%%%%%%%%%%%%%%%%%%%%%%%%%%%%%%%%%
% 2 %
\chapter{The pseudoholomorphic mapping construction \protect\\
on a higher jet level}
\markboth{2. Higher jet level}{2. Higher jet level}

\hspace{13.5pt}
In chapter~1 we found an algebra of invariants $\A^2_j=<j,N_j>$ generated
relative to commutator by the elements $j$, $N_j$, and such that
$\E^2_{x,y}=
\{\P:T_xL\to T_yM\ \vert \ \A^2_{j_M}\c\P^{\wedge\star}=\P\c\A^2_{j_L}\}$.
Applying the prolongation-projection method (see chapter~1) we contract the
equation $\E^2$ onto the equation $\E^\infty$ and find an algebra of
invariants $\A^\infty_j$ such that $\E^\infty_{x,y}=\{\P:T_xL\to T_yM\ \vert
\ \A^\infty_{j_M}\c\P^{\wedge\star}=\P\c\A^\infty_{j_L}\}$.
As it is shown in remark~7 (section 6.2) neither the algebra $\A^2_j$ nor
its closure
under differential Nijenhuis  bracket (see~\cite{FN}, \cite{ALV})
--- {\it the Nijenhuis algebra\/} $\A^\N_j$ --- coincides with the
algebra $\A^\infty_j$ looked for, i.e.
$\A^\N_j\ne\A^\infty_j=<j, N_j, \dots>$\,.
Below we construct an invariant coming from a higher jet space and being
a new (the third) generator of $\A^\infty_j$. Note that performing further
computations one could find other generators.

 \begin{th}\po
There exist new tensor invariants --- {\it higher Nijenhuis tensors\/}
--- such that $\A^\infty_j=<j, N_j, N_j^{(2)}, \dots>$\,.
The first of such invariants $N_j^{(2)}$ has the following form
\begin{equation}
 \begin{array}{l}
 \NN_j(\x,\e,\zeta,\nu)= -[N_j(\x,\e),jN_j(\zeta,\nu)]
-[jN_j(\x,\e),N_j(\zeta,\nu)] \phd                         \\
 +N_j([\x,jN_j(\zeta,\nu)],\e) +N_j(\x,[\e,jN_j(\zeta,\nu)])
 +jN_j([\x,N_j(\zeta,\nu)],\e) \phd                         \\
 +jN_j(\x,[\e,N_j(\zeta,\nu)])
 -N_j([\zeta,jN_j(\x,\e)],\nu) -N_j(\zeta,[\nu,jN_j(\x,\e)])
                               \phd                          \\
 -jN_j([\zeta,N_j(\x,\e)],\nu) -jN_j(\zeta,[\nu,N_j(\x,\e)]),\phd
% \tag{13}
\add
 \end{array}
\end{equation}
where $\x,\ \e,\ \zeta,\ \nu$ are arbitrary vector fields extending the
vectors prescribed at the point $x$. In general case this invariant
is the first among generating invariants from $\A_j^\infty$ distinguishing
structures $j$ and $(-j)$ (see remark~4).
 \end{th}

 \begin{proof}
As it is shown in theorem 3 for the equation
${\tilde\E}^1=\E^2\subset J^1_{PH}$ the projections ${\tilde\E}^k=\E^{k+1}$.
So instead of consideration the problem on projections of higher
prolongations of the equation $\E^1$, which leads to a formally integrable
equation (see chapter~1), one may consider the problem of finding out the
projections of higher prolongations of the new equation ${\tilde\E}^1$, and
at the first step, the problem of lifting on one level. The solution on
1-jet level of the equation ${\tilde\E}^1$ prolongs to the solution on
2-jet level iff there holds equality~(1)~and:
 \begin{eqnarray}
U(\x,\e,\te):& \! = & \!
\P^{(2)}(N_{j_L}(\x,\e),\te)-N_{j_M}(\P^{(2)}(\x,\te),\P\e)-
N_{j_M}(\P\x,\P^{(2)}(\e,\te))
\nonumber \\
{}& \! = & \! d(N_{j_M})(\P\x,\P\e,\P\te)- \P\c d(N_{j_L})(\x,\e,\te).
% \tag {14}
\add
 \end{eqnarray}
From (1) and (14) we have:
 $$
 \begin{array}{l}
U(\x,\e, j_L\te)+j_M U(\x,\e,\te)=
2j_M\c \P^{(2)}(N_{j_L}(\x,\e),\te) \phd                           \\
\quad +dj_M(\P\te,N_{j_M}(\P\x,\P\e)) -\P\c
dj_L(\te,N_{j_L}(\x,\e)) -N_{j_M}(dj_M(\P\te,\P\x),\P\e) \phd      \\
\quad -N_{j_M}(\P\x,dj_M(\P\te,\P\e))
+\P\c N_{j_L}(dj_L(\te,\x),\e)+\P\c N_{j_L}(\x,dj_L(\te,\e)) \phd  \\
\ =d(N_{j_M}) (\P\x,\P\e,j_M\P\te)+j_M\c d(N_{j_M})(\P\x,\P\e,\P\te)
\phd                                                               \\
\quad -\P\c d(N_{j_L})(\x,\e,j_L\te) -\P\c j_L\c d(N_{j_L})(\x,\e,\te).
\phd
 \end{array}
 $$
Hence:
 $$
2j_M\c\P^{(2)}(N_j(\x,\e),\te)= R_{j_M}(\P\x,\P\e,\P\te) -\P\c
R_{j_L}(\x,\e,\te),
 $$
where

 \begin{equation}
 \begin{array}{l}
R_j(\x,\e,\te)=  d(N_j)(\x,\e,j\te)+jd(N_j)(\x,\e,\te)        \phd  \\
\qquad\qquad
+N_j(dj(\te,\x),\e) +N_j(\x,dj(\te,\e)) -dj(\te,N_j(\x,\e)).  \phd
% \tag {15}
 \end{array}
\add
 \end{equation}

Since $\P^{(2)}$ is a symmetric tensor then for the tensor
 \begin{equation}
\NN_j(\x,\e,\zeta,\nu)=
R_j(\x,\e,N_j(\zeta,\nu)) -R_j(\zeta,\nu,N_j(\x,\e))
% \tag{16}
\add
 \end{equation}
we have:
 \begin{equation}
\NN_{j_M}\c(\P^\2)^\2=\P\c\NN_{j_L}.
% \tag {17}
\add
 \end{equation}

Therefore according to (17) we have a new invariant, which we will call
{\it a higher Nijenhuis tensor\/}, $\NN_j\in \A^3_j\subset\A^\infty_j$.
Formula (16) shows the
tensor nature of the invariant $\NN_j\in\wedge^2(\wedge^2T^*_x)\otimes T_x$,
while formula (13), following from formulae (14)-(15), shows its independence
of the Euclidean connection $\nabla$ which
defines $d=d_\nabla$. Note that a priori the tensor of the valence
(1,4) $\NN_j$ might become dependent with tensors of forth order coming from
$\A^2_j$ such as $N_j(N_j(\x,\e),N_j(\zeta,\nu))$, i.e. not to be new,
but example~2 shows that this is not the case. \hfill$\Box$
\end{proof}

 \begin{ex}{2}
Let $\p_i=\dfrac\p{\p x^i}$ be basis vectors on $\R^4=\R^4(x^1,x^2,x^3,x^4)$.
Let's define an almost complex structure by the formula:
 $$
j\p_1=\p_2,\ j\p_2=-\p_1,\ j\p_3=\p_4+x^2\p_1,\ j\p_4=-\p_3-x^2\p_2.
 $$
For the tensors $N_j$ and $\NN_j$ we have:
 $$
N_j(\p_1,\p_2)=0,\ \ N_j(\p_3,\p_4)=x^2\p_1,
 $$
 $$
N_j(\p_1,\p_3)=-N_j(\p_2,\p_4)=\p_1,\
N_j(\p_2,\p_3)=N_j(\p_1,\p_4)=-\p_2;
 $$
 $$
 \begin{array}{l}
\quad \NN_j(\p_1,\p_3,\p_3,\p_4)=-\p_1,
\phd                                                                 \\
\quad \NN_j(\p_1,\p_3,j\p_3,\p_4)= \NN_j(\p_1,\p_3,\p_4,\p_4)
+x^2\NN_j(\p_1,\p_3,\p_1,\p_4)=0.\phd\hfill {\Box}
 \end{array}
 $$
 \end{ex}

%%%%%%%%%%%%%%%%%%%%%%%%%%%%%%%%%%%%%%%%%%%%%%%%%%%%%%%%%%%%%%%%%%%%%%%%%%
% 3 %
 \chapter{A solvability criterion for the Cauchy-Riemann equation \protect\\
and formal normal forms of almost complex structures}
\markboth{}{}

%%%%%%%%%%%%%%%%%%%%%%%%%%%%%%%%%%%%%%%%%%%%%%%%%%%%%%%%%%%%%%%%%%%%%%%%%%
% 3.1 %
\section{Solvability \ of \  the \ Cauchy-Riemann \ %\protect\\
equation}
\markboth{3.1. Solvability of the Cauchy-Riemann equation}
{3.1. Solvability of the Cauchy-Riemann equation}

\hspace{13.5pt}
The set of \r-jets may be locally described as follows:

 \begin{eqnarray}
&& J^k_{PH}=\{(x,y,\P^{(1)},\d,\P^{(k)})\ | \ x\in L,y\in M, \quad
% \tag {18}
\add                                                          \\
&& \qquad\qquad\qquad\qquad\qquad\qquad\ \ \
\P^{(i)}\in Smbl^i_{PH}(L,M)\subset S^iT^*L\otimes TM\},\quad \nonumber
  \end{eqnarray}

 \begin{eqnarray*}
\sum\limits_
{\begin{array}{c}
\scriptstyle     r_1+\d+r_p=k \\
\scriptstyle     1\le p\le k
 \end{array}}
\!\!
\sum\limits_
{\begin{array}{c}
\scriptstyle     1\le i_s^r\le k, \#\{i_s^r\}=k \\
\scriptstyle     i_1^1<i_2^1<\d<i_p^1           \\
\scriptstyle     i_s^1<i_s^2<\d<i_s^{r_s}
 \end{array}}
\!\!\!\!\!
d^{p-1}j_M\left(\P^{(r_1)}(\x_{i_1^1},\d,\x_{i_1^{r_1}}),\d,
\P^{(r_p)}(\x_{i_p^1},\d,\x_{i_p^{r_p}})\right)=&&\\
\sum\limits_
{\begin{array}{c}
\scriptstyle     1\le p\le k, i_1^1=1, 1\le i_s^r\le k,\\
\scriptstyle     \#\{i_s^r\}=k, i_s^1<i_s^2<\d<i_s^*
 \end{array}}
\!\!\!\!\!\!
\P^{(k-p+1)}\left(dj_L^{p-1}(\x_{i_1^1},\d,\x_{i_1^p}),\x_{i_2^1},\d,
\x_{i_2^{k-p}}\right).&&
 \end{eqnarray*}
Here $\# S$ denotes the number of different elements of the set $S$.

The following theorem generalizes theorem 1 corresponding to the case $k=2$.

 \begin{th}\po
Let the mapping $u: \ol \to \om$ be such that
 \begin{equation}
j_M\c u_*=u_*\c j_L (\op{mod} \m^{k-1}),\
N_{j_M}\c u_*^\2\equiv u_*\c N_{j_L} (\op{mod} \m^{k-1})\
(k\ge2).
% \tag {19}
\add
 \end{equation}
Then there exists a mapping $\tilde u:\ol\to\om$ such that
 $$
j_M\c\tilde u_*-\tilde u_*\c j_L\in\m^{k},\ u_*-\tilde u_*\in\m^{k-1}.
 $$
In particular if $(u_*)_x\ne0$ then $({\tilde u}_*)_x\ne0$, and if $(u_*)_x$
is an embedding or a surjection then the same is true for $({\tilde u}_*)_x$.
 \end{th}

 \begin{proof}
A $(k-1)$-\r-jet
defined by means of symbols $(\P,\P^{(2)},\dots,\P^{(k-1)})$ is lifted
up to $k$-\r-jet iff some symmetric tensor $\P^{(k)}$ satisfies equation
(18). Consider the tensor
 \begin{equation}
P_k=j_M\c\P^{(k)}-\P^{(k)}\c(j_L\otimes\1_L^{k-1}).
% \tag{20}
\add
 \end{equation}

 \begin{lll}{1}
$(i_k)$\ \ \ $j_M\c P_k=-P_k\c(j_L\otimes\1_L^{k-1})$,
\newline
$(ii_k)$\ \ $P_k(\x,\e,\3)-P_k(\e,\x,\3)=
P_k(j_L\x,j_L\e,\3)-P_k(j_L\e,j_L\x,\3)$,
\newline
$(iii_k)$\  $P_k\in T^*_xL\otimes S^{k-1}T^*_xL\otimes T_yM$. \qed
 \end{lll}

 \begin{lll}{2}
Every tensor $P_k$, satisfying the
conditions of lemma $1_k$ may be represented in form (20), where
$\P^{(k)}\in S^kT^*_xL\otimes T_yM$.
 \end{lll}

 \begin{Proof}{Proof of lemma $2_k$}
Let's consider the tensor $B(\x,{\vec\e}^{(k-1)})$ such that $P_k=j_M\c
B-B\c(j_L\otimes \1_L^{k-1})$. Due to $(i_k)$ it always can be found,
and because of $(iii_k)$ one can suppose that $B$ is symmetric by the last
$(k-1)$ arguments: $B\in T^*_xL\otimes S^{k-1}T^*_xL\otimes T_yM$. Note that
under the change of $B$ by a tensor $j_L$-$j_M$-linear by the first argument,
the corresponding tensor $P_k$ does not change. Thus in the decomposition
$B=B^{0,*}+B^{1,*}$ by $j_L$-$j_M$-linear and antilinear by the first
argument components one may eliminate the first term and assume that
$B=B^{1,*}$.

Consider the decomposition
 $$
B^{1,*}=\sum\limits_{p=0}^{k-1}C_p,
 $$
where
$C_p=\sum\limits_{\tau_p\in\Z_2^{k-1}}B^{1,\tau_p}$, and $\Z_2=\{0,1\}$,
the number of units in $\tau_p$ is equal to $\#_1\tau_p=p$, and the tensor
$B^{1,\tau_p}$ is a $j_L$-$j_M$-linear or antilinear by $s$-th argument
if there's placed, respectively, 0 or 1 on the $s$-th place in $(1,\tau_p)$.
Note that if units form the set $i_1,\dots,i_r$ in $\tau_p$
and zeros form the set $j_1,\dots,j_{k-r-1}$, then the tensor $B^{1,\tau_p}$
is symmetric by the corresponding arguments $\x_{i_1},\dots,\x_{i_r}$ and by
$\x_{j_1},\dots,\x_{j_{k-r-1}}$. Note also that the tensor $B^{1,\tau_p}$
is symmetric by the arguments $\x_1,\x_{i_s}$. Actually, for example, for
$B^{1,1,*}$ from the formula $(ii_k)$ we have (we use
$j_L$-$j_M$-antilinearity by $\x_1$):
 $$
B^{1,*}(\x,\e,\3)+j_MB^{1,*}(\x,j_L\e,\3)=
B^{1,*}(\e,\x,\3)+j_MB^{1,*}(\e,j_L\x,\3),
 $$
 $$
\text{ i.e. }\quad
B^{1,1,*}(\x,\e,\3)=B^{1,1,*}(\e,\x,\3).
 $$
Thus the tensor $B^{1,\tau_p}$ is symmetric by arguments $\x_1,\x_{i_1},
\dots, \x_{i_r}$. It is not hard to show that if $\g\in S_{k-1}$ is a
permutation then the tensor $B^{1,\g(\tau_p)}$ is equal to
$B^{1,\tau_p}\c(\1\times\g)$, where $(\1\times\g)$ acts by permutation on
arguments; in the case when $\g(\tau_p)=\tau_p$ the statement coincides
with the noted above symmetry property of $B^{1,\tau_p}$. Hence if we add
the sum $\sum\limits_{\#_1\tau_{p+1}=p+1}B^{0,\tau_{p+1}}$ to the tensor
$C_p$, where $B^{0,\tau_{p+1}}=B^{1,\tau_p}\c\z$ for $\z\in S_k$ being such
a permutation that $(0,\tau_{p+1})=\z^{-1}(1,\tau_p)$, then we obtain a new
symmetric tensor ${\hat C}_p\in S^kT^*_xL\otimes T_yM$.

Let us consider the tensor $\P^{(k)}=\sum\limits_{p=0}^{k-1}{\hat C}_p$. It
is symmetric by the construction and $(\P^{(k)}-B)$ is $j_L$-$j_M$-linear.
Therefore the tensor $P_k=j_M\c\P^{(k)}-\P^{(k)}\c(j_L\otimes\1^{k-1})$
constructed by $\P^{(k)}$ is equal to the similar tensor constructed by
$B$, and lemma~$2_k$ is proved.

Note that the proposed construction slightly differs from the construction
of lemma~2 ($k=2$) in theorem~1. For the case $k=3$ we have:
 $$
B^{1,*}=B^{1,0,0}+(B^{1,0,1}+B^{1,1,0})+B^{1,1,1}.
 $$
The tensor $B^{1,1,1}$ is symmetric, and $B^{1,0,0}$ is symmetric by
the last two arguments $(\e,\te)$. Hence if $\z$ is a cyclic permutation
of arguments, $\z(\x,\e,\te)=(\e,\te,\x)$, then $B^{1,0,0}\c\z$ and
$B^{1,0,0}\c\z^2$ are $j_L$-$j_M$-linear by the first argument, and the
tensor $(B^{1,0,0}+B^{1,0,0}\c\z+B^{1,0,0}\c\z^2)$ is symmetric. The tensor
$B^{1,1,0}$ is symmetric by the arguments $(\x,\e)$ due to $(ii_3)$. So
$B^{1,1,0}(\x,\e,\te)=A(\x,\e,\te)+A(\e,\x,\te)$ for some tensor $A$
$j_L$-$j_M$-antilinear by the first two arguments and linear by the last.
Moreover, the tensor $B^{1,*}$ is invariant under the permutation $\tau$ of
the last two arguments, $\tau(\x,\e,\te)=(\x,\te,\e)$. Therefore
 \begin{eqnarray*}
&&B^{1,0,1}(\x,\e,\te)=B^{1,1,0}\c\tau(\x,\e,\te)=
A(\x,\te,\e)+A(\te,\x,\e), \qquad \\
&&B^{1,0,1}(\x,\e,\te)=B^{1,1,0}(\x,\te,\e)=
B^{1,1,0}(\te,\x,\e)=B^{1,1,0}\c\z^2(\x,\e,\te).\qquad
 \end{eqnarray*}
Let's define the tensor
 $$
B^{0,1,1}(\x,\e,\te):=B^{1,1,0}\c\z(\x,\e,\te)= A(\e,\te,\x)+A(\te,\e,\x).
 $$
Note that $(B^{1,1,0}+B^{1,0,1}+B^{0,1,1})$ is a symmetric tensor and
$B^{0,1,1}$ is $j_L$-$j_M$-linear by the first argument. So the tensor
 $$
P_k=j_M\c B^{1,*}-B^{1,*}\c(j_L\otimes\1\otimes\1)=
j_M\c\P^{(3)}-\P^{(3)}\c(j_L\otimes\1\otimes\1),
 $$
where $\P^{(3)}\in S^3T^*_xL\otimes T_yM$ is defined by the equality
 $$
\P^{(3)}=(B^{1,0,0}+B^{1,0,0}\c\z+B^{1,0,0}\c\z^2)+
(B^{1,1,0}+B^{1,0,1}+B^{0,1,1})+B^{1,1,1}.
 $$
\qed
 \end{Proof}

The statement of the lemmas may be reformulated as exactness of the
following complex:
 $$
0\to
{\cal G}_k \olra{i}
S^kT^*L\otimes TM \olra{\z_k}
T^*L\otimes S^{k-1}T^*L\otimes TM\olra{\k_k\oplus\nu_k}
(\Theta_k\oplus\Delta_k)\to0.
 $$
Here ${\cal G}_k$ is the prolongation of the symbol (\cite{ALV}, \cite{KLV})
of the Cauchy-Riemann equation $\E^1$, $S^kT^*L\otimes TM$ is the space of
all the symbols, $i$ is the inclusion, and $\z_k$ is the linearization of
the operator $\bar\partial$:
 $$
\z_k(\Psi)=j_M\c\Psi-\Psi\c(j_L\otimes\1^{k-1}).
 $$
The  maps $\k_k$ and $\nu_k$ on the spaces
 $$
\Theta_k=\wedge^2T^*L\otimes S^{k-2}T^*L\otimes TM \text{ and }
\Delta_k=S^{k-1}T^*L\otimes\op{Hom}_\C(TL,TM)
 $$
have the form:
 $$
\k_k=\delta_k\c{\tilde\z}_k \text{ and }
\nu_k(X\otimes{\vec Y}^{(k-1)}\otimes Z)=
({\vec Y}^{(k-1)}\otimes\epsilon(X\otimes Z)),
 $$
where $\delta_k$ is the operator of alternation by the first two arguments,
 $$
\delta_k(B)(X,Y,{\vec Z}^{(k-2)})=B(X,Y,{\vec Z}^{(k-2)})- B(Y,X,{\vec
Z}^{(k-2)}),
 $$
 $$
{\tilde\z}_k\Psi=j_M\c\Psi-\Psi\c(\1\otimes j_L\otimes\1^{k-2}),
 $$
and $\epsilon(B)=j_M\c B+B\c j_L$ is the linearization of the operator
$\partial$.
Note that the term $\Theta_k$, corresponding to the condition $(ii_k)$ of
lemma $1_k$, gives the preserving condition for the Nijenhuis tensor
differential and the term $\Delta_k$, corresponding to the condition $(i_k)$,
gives the identity. In more details, let us express the tensor $P_k$ by local
representative $(\P,\P^{(2)},\dots,\P^{(k-1)})$ of $(k-1)$-jet of the mapping
$u:\ol\to\om$ by means of formula (18) and let us substitute the
corresponding expression into the formulae of lemma $1_k$. Then the
conditions $(i_k)$ and $(iii_k)$ give identities and $(ii_k)$ gives the
preserving condition for $(k-1)$-jet of the Nijenhuis tensor.
The relevant calculations are quite similar to those of theorem~1
which correspond to the case $k=2$.

For example, if $k=3$ we have:
 \begin{eqnarray}
\!\!\!\!\!\!\!\!
&&\P_3(\x,\e,\te)=j_M\P^{(3)}(\x,\e,\te)-\P^{(3)}(j_L\x,\e,\te) \nonumber \\
\!\!\!\!\!
&&=\P d^2j_L(\x,\e,\te)-d^2j_M(\P\x,\P\e,\P\te)+
\P^{(2)}(dj_L(\x,\te),\e)+\P^{(2)}(dj_L(\x,\e),\te)             \nonumber \\
\!\!\!\!\!
&&-dj_M(\P^{(2)}(\x,\te),\P\e)-dj_M(\P\x,\P^{(2)}(\e,\te))-
dj_M(\P^{(2)}(\x,\e),\P\te).
% \tag{21}
\add
 \end{eqnarray}

The substitution of expression (21) into $(i_3)$ and $(iii_3)$ gives the
identity according to formulae (8) and the formula
 \begin{equation}
d^2j(j\x,\e,\te)=-jd^2j(\x,\e,\te)-dj(dj(\x,\te),\e)-dj(dj(\x,\e),\te).
% \tag{22}
\add
 \end{equation}

Let us substitute expression (21) into $(ii_3)$ and use formulae (1), (8),
(22) and the identity
 \begin{eqnarray*}
&&dj(N_j(\x,\e),\te)+jdN_j(\x,\e,\te)=\\
&&d^2j(j\x,j\e,\te)-d^2j(j\e,j\x,\te)-d^2j(\x,\e,\te)+d^2j(\e,\x,\te)+\\
&&dj(dj(\x,\te),j\e)-dj(dj(\e,\te),j\x)+dj(j\x,dj(\e,\te))-dj(j\e,dj(\x,\te)),
 \end{eqnarray*}
which is obtained from identity (9) by differentiation.
We have\footnote{We omit a page of calculations}:
\pagebreak

 \begin{eqnarray*}
\!\!\!\!\!\!\!\!\!\!\!
&&0=\PP(dj_L(\x,\te),\e)+\PP(dj_L(\x,\e),\te))-\\
\!\!\!\!\!\!
&&dj_M(\PP(\x,\te),\P\e)-dj_M(\P\x,\PP(\e,\te))-dj_M(\PP(\x,\e),\P\te)+\\
\!\!\!\!\!\!
&&\P d^2j_L(\x,\e,\te)-d^2j_M(\P\x,\P\e,\P\te)-\\
\!\!\!\!\!\!
&&\PP(dj_L(\e,\te),\x)-\PP(dj_L(\e,\x),\te)+\\
\!\!\!\!\!\!
&&dj_M(\PP(\e,\te),\P\x)+dj_M(\P\e,\PP(\x,\te))+dj_M(\PP(\e,\x),\P\te)-\\
\!\!\!\!\!\!
&&\P d^2j_L(\e,\x,\te)+d^2j_M(\P\e,\P\x,\P\te)-\\
\!\!\!\!\!\!
&&\PP(dj_L(j_L\x,\te),j_L\e)-\PP(dj_L(j_L\x,j_L\e),\te)+\\
\!\!\!\!\!\!
&&dj_M(\PP(j_L\x,\te),j_M\P\e)+dj_M(j_M\P\x,\PP(j_L\e,\te))+
 dj_M(\PP(j_L\x,j_L\e),\P\te)\\
\!\!\!\!\!\!
&&-\P d^2j_L(j_L\x,j_L\e,\te)+d^2j_M(j_M\P\x,j_M\P\e,\P\te)+\\
\!\!\!\!\!\!
&&\PP(dj_L(j_L\e,\te),j_L\x)+\PP(dj_L(j_L\e,j_L\x),\te)-\\
\!\!\!\!\!\!
&&dj_M(\PP(j_L\e,\te),j_M\P\x)-dj_M(j_M\P\e,\PP(j_L\x,\te))
-dj_M(\PP(j_L\e,j_L\x),\P\te)\\
\!\!\!\!\!\!
&&+\P d^2j_L(j_L\e,j_L\x,\te)-d^2j_M(j_M\P\e,j_M\P\x,\P\te)\\
\!\!\!\!\!\!\!\!\!\!
&&=j_M[-\PP(N_{j_L}(\x,\e),\te)-\P dN_{j_L}(\x,\e,\te)+\\
\!\!\!\!\!\!
&&N_{j_M}(\PP(\x,\te),\P\e)+N_{j_M}(\P\x,\PP(\e,\te))+
d N_{j_M}(\P\x,\P\e,\P\te)].
 \end{eqnarray*}
But the vanishing of the last expression means exactly that the 1-jet
from ${\tilde\E}^1$ may be prolonged to a 2-jet.

We omit calculations for the case of general $k$. The theorem is proved.
\hfill
$\Box$
 \end{proof}

 \begin{cor}
Let the Nijenhuis tensors be of smallness of
order $k$ at the point $x\in L$ and $y\in M$:
$N_{j_L}\in\m^k_x$, $N_{j_M}\in\m^k_y$ ($1\le k\le\infty$). Then there
exists such a mapping $u: \ol\to\om$ that $du_x\ne0$ and
 $$
\hfill\hfill\qquad\quad\quad\quad\quad\quad\quad\quad\quad\quad\
j_M\c u_*-u_*\c j_L\in\m^{k+1}.
\quad\quad\qquad\quad\quad\quad\quad\quad\quad\quad\ \hfill\hfill
\Box
 $$
 \end{cor}

 \begin{rk}{2}
When $\op{dim}_\C L=\op{dim}_\C M=n$, $k=\infty$,
this assertion is the formal part of Newlander-Nirenberg theorem
(see~\cite{NN}, \cite{NW}), stating that if the Nijenhuis tensors of the
almost complex structures $j$ and $j_0$ coincide then these structures
are locally isomorphic; here $j_0$ is the standard complex structure
in $\ok_{\C^n}(0)$ (certainly $N_{j_0}\equiv0$, see formula (9)). A
generalization of this fact is theorem~6 below. Note also that smooth
equivalence does not follows from the formal one. Actually, one may
construct an almost complex structure $j$ such that $N_j\in\m^\infty$,
but $N_j\not\equiv0$ (see the proof of theorem~5 at chapter~4).
Then $j$ is formally equivalent to $j_0$ but not equivalent smoothly. \qed
 \end{rk}

 \begin{rk}{3}
The conditions of the corollary might be reformulated in the following form:
$dN_{j_L}\in\mu_x^{k-1}$, $dN_{j_M}\in\mu_y^{k-1}$, where
$d=d_\nabla$ and $\nabla$ is a trivial connection. Let us note that
there are no almost complex structures save complex ones such that
$d(N_j)\equiv0$. Actually, in this case there exist coordinates in which
the components of the Nijenhuis tensors are constant, and additionally we
may suppose that the induced coordinates in the fixed tangent space
$T_x\simeq\C^n$ are complex. Let's write the action of the almost complex
structure on the basis vectors $\p_r=\frac\p{\p x^r}$:
 $$
j\p_{2r-1}=\p_{2r}+\sum a^i_{2r}\p_i,\
j\p_{2r}  =-\p_{2r-1}+\sum a^i_{2r-1}\p_i;\ \
a^s_t\in\m_x.
 $$
If for almost any vector $\xi$: $\op{dim}N_j(\x,T_x)=2n-2$, i.e. the
tensor $N_j$ is of general position, ${\n}_{N_j}=2$ almost everywhere,
in the sense of definition~4 from chapter~4, then from
$N_j(j\p_{2r-1},\p_{2r})\equiv N_j(j\p_{2s-1},\p_{2s})\equiv0$ it follows
that $a_{2s}^{2r}\equiv a_{2s-1}^{2r}\equiv a_{2s}^{2r-1}\equiv
a_{2s-1}^{2r-1}\equiv 0$ for all $r\ne s$, i.e. by definition
$N_j(\p_{2r},\p_{2s})=0$; contradiction. Let us suppose now that there is a
subspace $K_1\subset T_x$ ($\op{dim}K_1\ge4$) such that
$\left.N_j\right\vert_{\wedge^2K_1}\equiv0$, and we also suppose that $K_1$
is maximal satisfying this property. Extending $K_1$ by affine shifts
in the fixed coordinate system to a neighborhood of the point $x$ we obtain
a foliation ${\cal K}_1$, $\op{dim}{\cal K}_1=\op{dim}K_1$, $j$-invariant
according to $N_j(j\x,\e)=-jN_j(\x,\e)=0$ for $\x, \e\in T_\bullet{\cal K}_1$,
and the restrictions of $j$ to fibers of which are integrable,
$\left.N_j\right\vert_{{\cal K}_1}\equiv0$. In an additional subspace
let's choose a maximal degenerate relative to $N_j$ subspace $K_2$, and
construct the foliation ${\cal K}_2$ by it and so on.
Finally, we are lead to the situation when the almost complex structure
$j$ in a neighborhood $\ok(x)$ is represented as a direct sum of complex
structures on the fibers of the standard foliations ${\cal K}_i$,
which are obtained from the summand of the complex decomposition
$T_x=\oplus K_i$. Now it is easy to obtain contradiction with the condition
$dN_j\equiv0$, if $j$ is not globally integrable.

Let us note however that there are nontrivial connections $\nabla$ such that
$d_\nabla(N_j)\equiv0$. Actually, if one takes for $\nabla$ an almost
complex connection (\cite{L}), then $d_\nabla(j)\equiv0$, which implies
$d_\nabla(N_j)\equiv0$. In particular, this property is satisfied by
a minimal connection $\nabla$, which is nonsymmetric and have the torsion
tensor equal to $T_\nabla=N_j$; such connections always exist, see~\cite{L}.
\qed
 \end{rk}

%%%%%%%%%%%%%%%%%%%%%%%%%%%%%%%%%%%%%%%%%%%%%%%%%%%%%%%%%%%%%%%%%%%%%%%%%%
% 3.2 %
\section{Formal normal forms %\protect\\
of almost complex structures}
\markboth{3.2. Formal normal forms}{3.2. Formal normal forms}

\hspace{13.5pt}
By means of theorem 3 it is possible to construct
{\it formal normal forms of almost complex structures\/}.
One may consider almost complex structures $modulo\,\mu^k$, i.e. such
automorphisms $j: T_{x^\prime}X\to T_{x^\prime}X$,
$x^\prime\in X$, that $j^2=-\1_X (\op{mod}\,\mu^k_x)$.

Let us suppose that we constructed normal forms of almost complex structures
$modulo\,\mu^k$: $\J_k=\{j_{\ll_1,\d,\ll_{k-1}}\}$. Here
$\ll_1,\d,\ll_{k-1}$ are parameters determining elements $\J_k$. Consider
those of the structures that prolong to almost complex structures
$modulo\,\mu^{k+1}$, i.e. those $J_k\in\J_k$ that there exists an
automorphism of the tangent bundle $\D J_k\in\mu^k$ for which
the structure $J_k+\D J_k$ is almost complex $modulo\,\mu^{k+1}$; let us
call such structures $k$-compatible. If $J_k$ is almost complex structure
$modulo\,\mu^{k+1}$ and $\D J_k\in\mu^k$, then $J_k+\D J_k$ is almost
complex structure $modulo\,\mu^{k+1}$ iff
$j_0\c\D J_k+\D J_k\c j_0\in\mu^{k+1}$,
where $j_0=J_k(\op{mod}\mu)$ is the linearization of the structure $J_k$.

 \begin{prop}\po
Two almost complex $modulo\,\mu^{k+1}$ ($k\ge1$)
structures $J_k$ and $J_k+\D J_k$ are equivalent by module $\mu^{k+1}$ iff
the automorphism $\D J_k$ is compatible with the complex structure
$j_0$ by module $\mu^{k+1}$:
 \begin{eqnarray*}
N_{(j_0,\D J_k)}(\x,\e)
& \!\!\! := \!\!\! &
[j_0\x, \D J_k\e] + [\D J_k\x,j_0\e]
- j_0[\x, \D J_k\e]                                 \\
& \!\!\! - \!\!\! &
j_0[\D J_k\x, \e] - \D J_k[\x,j_0\e] - \D J_k[j_0\x,\e]\in\mu^{k+1}.
 \end{eqnarray*}
 \end{prop}

Actually, $N_{(J_k+\D J_k,J_k)}:= N_{(J_k+\D J_k)}- N_{J_k}\equiv
N_{(j_0,\D J_k)}(\op{mod}\,\mu^{k+1})$.\hfill $\Box$

Thus, $k$-compatible almost complex structures
$j_{\ll_1,\d,\ll_{k-1}}$ from $\J_k$ define normal forms of almost complex
structures $modulo\,\mu^{k+1}$ $j_{\ll_1,\d,\ll_k}$ from $\J_{k+1}$, where
the parameter $\ll_k$ takes values in the symbols of $k$-jets of the
compatibility Nijenhuis tensors
 $$
\ll_k= [N_{(j_0,\D J_k)}]_k\in
(\mu^k\wedge^2T^*_x\otimes T_x)/(\mu^{k+1}\wedge^2T^*_x\otimes T_x)
=\wedge^2T^*_x\otimes S^kT^*_x\otimes T_x.
 $$
(compare the values space with the term $\Theta_{k+2}$ in the complex from
the proof of theorem~3).

So we have constructed normal forms of almost complex structures\linebreak
$modulo\,\mu^{k+1}$: $\J_{k+1}=\{j_{\ll_1,\d,\ll_k}\}$, i.e. we have found
such structures that for every almost complex structure ($modulo\,\mu^{k+1}$)
$j:\ok(x)\to\ok(x)$ there exists a diffeomorphism $\vp_k:\ok(x)\to\ok(x)$
and a collection of parameters $\ll_1,\d,\ll_k$, which satisfy the identity
 $$
\vp^*_kj\equiv j_{\ll_1,\d,\ll_k}(\op{mod}\,\mu^{k+1}).
 $$

Let's call {\it the tree of almost complex structures\/} the hierarchy
$\J_1=\{j_0\}$, $\J_2=\{j_{\ll_1}\}$, $\J_3=\{j_{\ll_1,\ll_2}\}$, and so on.
Let's call {\it a branch of the tree\/} any sequence of prolonging each
other almost complex structures from the tree $j_0$, $j_{\ll_1}$,
$j_{\ll_1,\ll_2}$, $j_{\ll_1,\ll_2,\ll_3}$, $\d$; let us denote such a
sequence by the symbol $j_{\{\ll_i\}_1^\infty}$, and their union by the
symbol $\J_\infty$. By Borel's theorem for any formal diffeomorphism there
exists a smooth one inducing it. Hence from the arguments above it follows
the statement:

 \begin{th}\po
The set $\J_\infty$ defines formal normal forms of almost complex structures
$modulo\,\mu^\infty$, i.e. for any formal almost complex structure
$j$, $j^2+\1\in\mu^\infty$, there exists such an element
$j_{\{\ll_i\}_1^\infty}\in\J_\infty$ and a diffeomorphism
$\vp:\ok(x)\to\ok(x)$, that
 $$
\hfill\hfill\quad\qquad\quad\quad\quad\quad\quad\quad\quad\ \
\vp^*j\equiv j_{\{\ll_i\}_1^\infty}(\op{mod}\,\mu^\infty).
\quad\quad\quad\qquad\quad\quad\quad\quad\quad\ \ \hfill\hfill
\Box
 $$
 \end{th}

Note that the constructed normal forms are not minimal in the sense of
decomposition on nonintersecting classes even $modulo\,\mu^2$, i.e. by
module of Nijenhuis tensors classification at a point, see the arguments
after theorem~1. But nevertheless this classification essentially reduces
the number of possible Taylor series expansions of almost complex
structures, e.g. there is only one branch of the tree which correspond to
all complex structures.

%%%%%%%%%%%%%%%%%%%%%%%%%%%%%%%%%%%%%%%%%%%%%%%%%%%%%%%%%%%%%%%%%%%%%%%%%%
% 4 %
 \chapter{Linear theory of Nijenhuis tensors}
\markboth{4. Linear theory of Nijenhuis tensors}
{4. Linear theory of Nijenhuis tensors}

\hspace{13.5pt}
Let us consider a Nijenhuis tensor $N=N_j$ on a fixed tangent space
$T_x$ as a bilinear map $N: T_x\otimes T_x\to T_x$. This map is antilinear
and skew-symmetric, $N\in\N_j$, see (10). There're no other constraints on
it.

 \begin{th}\po
For any tensor $N\in\N_{j_0}(T_x)$ (see (10))
there exists an almost complex structure $j$ in $\ok(x)$ such that
$j_x=j_0$ and $N=(N_j)_x$.
 \end{th}

 \begin{proof}
Let $(x^i)$ be such a coordinate system in a neighborhood of the point $x$
that the vectors $\p_i=\left.\dfrac\p{\p x^i}\right\vert_x$ form the standard
complex relative to $j_0$ basis in $T_x$. Then in a neighborhood of $x\in X$:
 \begin{equation}
  \begin{array}{l}
j\p_{2r-1}=\p_{2r}+\sum a^i_{2r}\p_i                      \phd     \\
\quad
j\p_{2r}  =-\p_{2r-1}+\sum a^i_{2r-1}\p_i;\ \ a^i_s(0)=0; \phd
% \tag {23}
\add
  \end{array}
 \end{equation}
we suppose the coordinates $x^i=0$ at the point $x$. We have:

 \begin{eqnarray}
\!\!\!\!N_j(\p_{2s-1},\p_{2t-1})_0
& \!\!\!\!\! = \!\!\!\!\! &
\sum\left( \frac{\p a_{2t}^{2i-1}}{\p x^{2s}}(0)
-\frac{\p a_{2s}^{2i-1}}{\p x^{2t}}(0)
+\frac{\p a_{2t}^{2i}}{\p x^{2s-1}}(0)
-\frac{\p a_{2s}^{2i}}{\p x^{2t-1}}(0)\right)\!\p_{2i-1}  \nonumber  \\
& \!\!\!\!\! + \!\!\!\!\! &
\sum\left( \frac{\p a_{2t}^{2i}}{\p x^{2s}}(0)
-\frac{\p a_{2s}^{2i}}{\p x^{2t}}(0)
-\frac{\p a_{2t}^{2i-1}}{\p x^{2s-1}}(0)
+\frac{\p a_{2s}^{2i-1}}{\p x^{2t-1}}(0)\right)\p_{2i}  \nonumber  \\
& \!\!\!\!\! = \!\!\!\!\! &
\sum\left( c_{s,t}^{2i-1}\p_{2i-1}+c_{s,t}^{2i}\p_{2i} \right).
% \tag {24}
\add
 \end{eqnarray}

The theorem states that under arbitrary choice of the constants $c^i_{s,t}$,
$1\le i\le 2n$, $1\le s<t\le n$, equation (24) has a solution. Actually,
let $a^i_{2t}(x^1,\d,x^{2n})=\sum\limits_{s<t} c_{s,t}^i x^{2s}$ and let's
define $a^i_{2t-1}$ from the condition $j^2=-\1$, using formula
(23).\ \hfill $\Box$
 \end{proof}

 \begin{dfn}{3}
Let's call {\it the space of linear Nijenhuis tensors\/} the set of tensors
$\N_j^{(n)}$ on a fixed complex linear space $V\simeq \C^n$, defined by
formula (10).
 \end{dfn}

Let us consider the map associating to a Nijenhuis tensor $N\in\N_j(V)$ on
$V\simeq\C^n$ the function $\n_N: V\to2\Z_+$:
 $$
\n_N(\x):=\op{dim}\{\e\ \vert\ N(\x,\e)=0\}=
2\op{dim}_\C\{\e\ \vert\ N(\x,\e)=0\}\in 2\Z_+.
 $$
Since $\n_N(\x)=\n_N(j\x)$ we could also consider the function $\n_N$ as
the mapping $\n_N: PV\simeq \C P^{n-1}\to 2\Z_+$.
The function $\n_N(\x)$ is upper semicontinuous by $\x$:
$\overline{\lim\limits_{\x\to\x_0}}
\n_N(\x)\le \n_N(\x_0),\ \forall \x_0\in V$;
it's also upper semicontinuous by $N$:
$\overline{\lim\limits_{N^\prime\to N}}
\n_{N^\prime}(\x)\le \n_N(\x),\ \forall N\in\N$.

 \begin{dfn}{4}
Let us say that a linear Nijenhuis tensor $N\in\N_j$ is of general position
if for almost all $\x$ we have $\n_N(\x)=2$, i.e. restrictions of $N$ onto
two-dimensional planes vanish only on complex lines save the set of
two-dimensional planes of zero measure.
 \end{dfn}

It is easy to show that the condition in definition~4 is really a
condition of general position: the set of tensors $N\in\N_j$ with
$\n_N(\x)=2$ for a.e. $\x$ is open due to semicontinuity and everywhere
density is easily proved by using example~3 below (this is proved in
details in appendix).

 \begin{ex}{3}
Let us consider a bilinear skew-symmetric map
$A:\wedge^2(\R^n)\to\R^n$, defined by the formula:
 $$
A(\x,\e)=\sum_{i=1}^n(\x^i\e^{i+1}-\x^{i+1}\e^i)e_i,\,
\x=\sum_{i=1}^n\x^ie_i,\, \e=\sum_{i=1}^n\e^ie_i,\, \x^{n+1}=\x^1,\,
\e^{n+1}=\e^1.
 $$
The tensor $A$ satisfies the property that for a.e. vector $\x$ it follows
from the condition $A(\x,\e)=0$ that vectors $\x$ and $\e$ are parallel.

Define a complex structure $j$ on $\R^{2n}=\R^n_1\oplus\R^n_2$ by the
formula $j(\x\oplus\e)=(-\e,\x)$, and define the corresponding
antiinvariant tensor $N\in\N(\R^{2n})$, $N:\wedge^2(\R^{2n})\to\R^{2n}$
by the formula:
 $$
N(\x_1\oplus\e_1,\x_2\oplus\e_2)= [A(\x_1,\x_2)-A(\e_1,\e_2)]\oplus
[-A(\x_1,\e_2)-A(\e_1,\x_2)].
 $$
It's easy to see that $\n_N(\x)=2$ for a.e. $\x$ (see appendix for details).
\hfill\qed
 \end{ex}

It follows from definition 4 that complex lines might be singled out from
two-dimensional planes by means of a Nijenhuis tensor of general position
$N$: $\C\x=\op{Ker}N(\x,\cdot)$ for almost all $\x\ne0$, and for others use
continuity.
Hence one might define complex independency of vectors using a Nijenhuis
tensor $N$ of general position: vectors $\x_1,\d,\x_r$ are complex
independent if for any vector $\x$ from the linear span
$\Lambda=<\x_1,\d,\x_r>$
and any number $s$ the equality $N(\x_s,\x)=0$ implies
$\x=\op{const}\cdot\x_s$, or the same for a small perturbation of $\Lambda$.
Moreover, from any complex independent collection of $n=\op{dim}_\C V$
vectors in $V$ one can canonically construct by means of the tensor $N$
a decomposition on complex lines $V=\oplus_{i=1}^n\C$, associating to each
vector its (two-dimensional) $N$-annihilator (or doing the same after a small
perturbation and then taking the limit when the perturbation vanishes).

 \begin{prop}
Let $N\in\wedge^2 V^*\otimes V$ be a skew-symmetric
tensor, antilinear regardness two almost complex structures $j_1$ and $j_2$
on a complex linear space $V$, $N\in \N_{j_1}\cap\N_{j_2}$. If $N$ is a
tensor of general position then $j_1=\pm j_2$.
 \end{prop}

 \begin{proof}
Let us consider arbitrary complex line $\C\subset V$, i.e. two-dimensional
complex space, the restriction of the Nijenhuis tensor $N$ on which
vanishes, if we suppose that it is complex generated by the vector $\x$ such
that $\n_N(\x)=2$, this is the case of general position due to definition~4;
the invariance of this space regarding to $j_1$ and $j_2$ follows
automatically. There exists a vector $\x\ne0$ on this line such that
the vectors $j_1\x$ and $j_2\x$ are parallel. Actually, up to dilations,
the complex multiplication $j_k$ acts on a complex line by rotations:
 $$
\vp\mapsto\vp+\delta_k(\vp)\mapsto\vp+\pi,\ \ \vp\in S^1=\R^1(\op{mod}\,2\pi).
 $$
Changing $j_2$ by $(-j_2)$, if needed, one may suppose that
$\delta_1(\vp),\delta_2(\vp)\in(0,\pi)$. Neither of inequalities
$\delta_1(\vp)<\delta_2(\vp)$, $\delta_1(\vp)>\delta_2(\vp)$ could be true
for all $\vp$ since $j_1^2=j_2^2=-\1$. Therefore, there exists $\vp$ such
that $\delta_1(\vp)=\delta_2(\vp)$, i.e. there exists a vector $\x\ne0$
for which $j_1\x=\a j_2\x$, $\a\in\R\setminus\{0\}$.

For any vector $\e\notin\C=\C\x$ we have:
 $$
N(\x,j_1\e)=N(j_1\x,\e)=N(\a j_2\x,\e)=N(\x,\a j_2\e),
 $$
i.e. $N(\x,j_1\e-\a j_2\e)=0$. But $j_1\e\in\C\e$ and $j_2\e\in\C\e$ (see
the remarks before the proposition that $\C\e$ depends not on the
structures $j_1$ and $j_2$), and also $\C\x\cap\C\e=\{0\}$.
So due to the condition that $\n_N(\x)=2$ almost everywhere we have $j_1\e=\a
j_2\e$ for any vector $\e\notin\C\x$, and hence, by continuity, also for any
$\e\in V$. Thus, $j_1=\a j_2$, from where $\a=\pm1$. \hfill $\Box$
 \end{proof}

Let us consider a linear space $V$, two complex structures $j_1$ and $j_2$
on it and arbitrary Nijenhuis tensor $N\in\N_{j_1}\cap\N_{j_2}$, not
necessary of general position. Suppose that $N\not\equiv0$. Then there
exists a vector $\te=N(\x,\e)\ne0$. Note that
$j_1j_2\te=N(j_1\x,j_2\e)=j_2j_1\te$. Therefore the subspace
$\T= {<\te, j_1\te, j_2\te, j_1j_2\te>}$ is $j_1$- and $j_2$-invariant.
Note that $\T$ is decomposed into the sum of invariant subspaces:
 $$
\T=<\te^+, j_1\te^+> \oplus <\te^-,j_1\te^->,
 $$
where $\te^+ =j_1\te +j_2\te$
and $\te^-=j_1\te-j_2\te$ satisfy the properties $j_1\te^+ =j_2\te^+$,
$j_1\te^-= -j_2\te^-$. Hence, if the subspace $\T$ is four-dimensional then
there exist nonzero vectors $\te^\pm$ such that $j_1\te^\pm=\pm j_2\te^\pm$;
if on the other side $\T$ is two-dimensional then there exists only one of
such vectors. Let us denote
 $$
\Pi^+=\{\x\,\vert\,j_1\x=j_2\x\},\ \Pi^- =\{\x\,\vert\, j_1\x=-j_2\x\},\
\Pi=\Pi^+\oplus\Pi^-.
 $$
We showed that any vector
$\te\in\op{Im}N(\cdot,\cdot)$ is represented as a sum of a vector from
$\Pi^+$ and a vector from $\Pi^-$. Note that any vector $\te\in V$ is
represented as a sum $\te= \te_{(+)}+\te_{(-)}$, where
$\te_{(\pm)} =\frac12(\te \mp j_1j_2\te)$, where $N(\te_{(+)},\Pi^-)=0$,
$N(\te_{(-)},\Pi^+)=0$. Let us denote
 $$
K^+= <\x\,\vert\,N(\x,\Pi^-)=0>,\
K^-= <\x\,\vert\,N(\x,\Pi^+)=0>.
 $$
We have proved:

 \begin{prop}
For the subspace $\hat\Pi \subset V$, generated by
$\op{Im}N(\cdot,\cdot)$, the following inclusion holds: $\hat\Pi\subset\Pi$.
Moreover $\Pi^+\subset K^+$, $\Pi^-\subset K^-$, $V=K^++K^-$.
The intersection $K^+\cap K^-$ coincides with the space of vectors
$K=\op{Ker}N(\cdot,\Pi) :=\{\x\,\vert\,N(\x,\Pi)=0\}$.\qed
 \end{prop}

 \begin{ex}{4}
Let us consider a linear space $V^8=\C_1\oplus\C_2\oplus\C_3\oplus\C_4$,
$\C_s=$ $<\p_{2s-1},\p_{2s}>$, with the standard complex structure $j_1$ and
with a structure $j_2$, which differs from $j_1$ on basis vectors only on
$\C_3$, where $j_2\p_5=\p_6+\p_1$, $j_2\p_6=-\p_5-\p_2$. Let's consider the
tensor $N\in\N_{j_1}\cap\N_{j_2}$, given by the conditions:
 $$
N(\p_1,\p_3)=\p_1,\ N(\p_5,\p_7)=\p_1,\ N(\C_1\oplus\C_2, \C_3\oplus\C_4)=0.
 $$
We have: $\Pi^-=0$, $\Pi=\Pi^+=\C_1\oplus\C_2\oplus\C_4$,
$K^-=\op{Ker}N(\cdot,\Pi)=\C_4$, $K^+=V^8$, $\op{Ker}N(\cdot,\cdot)=0$.
\qed
 \end{ex}

%%%%%%%%%%%%%%%%%%%%%%%%%%%%%%%%%%%%%%%%%%%%%%%%%%%%%%%%%%%%%%%%%%%%%%%%%%
% 5 %
 \chapter{Almost complex structures of general position \protect\\
and their classification}
\markboth{5. Almost complex structures of general position}
{5. Almost complex structures of general position}

\hspace{13.5pt}
In this chapter we study the general properties of Nijenhuis tensors and
simplify the corollary from theorem~3 (chapter~3)
for almost complex structures of general position.

 \begin{dfn}{5}
We say that an almost complex structure $j$ on $\ok(x)$ is of {\it general
position\/} at a point $x$, if for the jet of the Nijenhuis tensor
the following separability property holds: for almost every vector
$\x\in T_x$ there exists a number $k=k(\x)\ge0$ such that with some
trivial connection $\nabla$ it follows from
$d^k_\nabla N_j(\x,\e,{\vec\nu}^{(k)})=0$, $\e\in T_x$, for all
${\vec\nu}^{(k)}\in S^kT_x$, that $\e\in\C\x=<\x,j\x>$.
Moreover we suppose that the pair $(\ok(x),j)$ is not isomorphic to the pair
$(\ok(x),-j)$.
 \end{dfn}

Definition 4 for linear Nijenhuis tensor is obtained from definition~5 when
$k\equiv0$. It's easy to see that the set of almost complex structures
$j$ on $\ok(x)$ of general position contains an open and everywhere dense
set in the set of all almost complex structures; hence it follows that the
set of the points, where an almost complex structure is of general position,
is open. Note that similar to the linear case for an almost complex
structure $j$ of general position from any complex independent
collection of $n$ vectors in $T_x$ one can canonically construct by means of
the tensor $N_j$ a complex polarization: $T_x=\oplus_{i=1}^n\C$.

 \begin{rk}{4}
An almost complex structure $j$ of general position on $\ok(x)$ is not
isomorphic to the structure $(-j)$. Actually, for the invariant, introduced
in chapter~2, we have: $\NN_{-j}=-\NN_j$, and for a general almost complex structure
$j$ tensors $\NN_j$ and $(-\NN_j)$ (on the fixed tangent space) are
nonisomorphic. Note that $(T_x,j_x)\simeq(T_x,-j_x)$, and $N_{-j}=N_j$,
from where it follows that $j \simeq (-j)\,(\op{mod}\,\mu^2_x)$ on $\ok(x)$,
but almost always $j \not\simeq (-j)\,(\op{mod}\,\mu^3_x)$ and the
tensor $\NN_j$ is the first generating invariant in $\A_j^\infty$,
which distinguish, in general case, almost complex structures $j$ and $(-j)$.
\qed
 \end{rk}

 \begin{ex}{5}
Consider the almost complex structure in $\R^4$,
defined by the equalities:

 \begin{eqnarray*}
j\p_1=\p_2,{\ \ } && \qquad j\p_3=\p_4+(x^2)^2\p_1+\ve(x^3)^2\p_2, \\
j\p_2 =-\p_1,     && \qquad j\p_4=-\p_3-(x^2)^2\p_2+\ve (x^3)^2 \p_1.
 \end{eqnarray*}
For $\ve=0$ almost complex structures $j$ and $(-j)$ are isomorphic, and
for $\ve\ne0$ are not. It's easy to see that the given almost complex
structure for any $\ve\ne0$ is of general position.\qed
 \end{ex}

 \begin{th}\po
Let us suppose that almost complex structures $j_L$ on
$\ol$ and $j_M$ on $\om$ are of general position.

$1^\c$.
Let us suppose that there exists a mapping $u:\ol\to\om$, which conjugates
$N_{j_L}$ with $N_{j_M}$ and is a diffeomorphism on the image.
Then $u$ is either pseudoholomorphic or antipseudoholomorphic mapping:
$j_M\c u_*=\pm u_*\c j_L$.

$2^\c$.
If the Nijenhuis tensors are conjugated formally,
$N_{j_M}\c u_*^\2 \equiv u_*\c N_{j_L}$ \linebreak
$(\op{mod}\m^\infty_{x,y})$, and
the dimensions of manifolds $L$ and $M$ are equal, $2l=2m$, then the almost
complex structure $j_L$ is formally equivalent to exactly one of the
structures $j_M$ or $(-j_M)$.
 \end{th}

 \begin{proof}
Let us consider the case of smooth conjugacy of the Nijenhuis tensors; the
formal one is obtained by passing on to series. Let
$N_{j_M}\c\wedge^2u_*= u_*\c N_{j_L}$.
One may choose complex parallelizations
$T_{x^\prime}L=\oplus_1^l\C$ and $T_{y^\prime}M=\oplus_1^m\C$ such that
the mapping $u_*$ includes one into the other. Therefore the image $u(L)$
is $j_M$-invariant. Hence one can consider the case of the manifolds of the
same dimensions $2l=2m$ and suppose that $u$ is a diffeomorphism. Thus we
have at a point $y^\prime\in u(L)\cap M$ the Nijenhuis tensor
$N_{j_M}=u_*\c N_{j_L}\c (u_*^\2)^{-1}$ and two complex structures $j_M$ and
$u_*\c j_L\c u_*^{-1}$ on $T_{y^\prime}$, regardness to which it is
antilinear. Note that the statement of proposition~2 holds true if instead
of the tensor $N$ one considers the jet $\oplus_k d^kN_j$ of the Nijenhuis
tensor of general position. Therefore, $j_M=\pm u_*j_L$, and since an
almost complex structure $j$ of general position is not isomorphic to
$(-j)$, the theorem is proved. \hfill $\Box$
 \end{proof}

 \begin{rk}{5}
The theorem is not true for almost complex structures not of general
position. If $j_L$ on $\ol$ and $j_M$ on $\om$ are structures of general
position then the Nijenhuis tensors of almost complex structures
$j_L\oplus j_M$ and $j_L\oplus(-j_M)$ on $\ol\t\om$ coincide, while
the structures themselves are neither isomorphic nor antiisomorphic.
\qed
 \end{rk}

%%%%%%%%%%%%%%%%%%%%%%%%%%%%%%%%%%%%%%%%%%%%%%%%%%%%%%%%%%%%%%%%%%%%%%%%%%
% 6 %
 \chapter{Two special cases}
\markboth{6. Two special cases}{6. Two special cases}

%%%%%%%%%%%%%%%%%%%%%%%%%%%%%%%%%%%%%%%%%%%%%%%%%%%%%%%%%%%%%%%%%%%%%%%%%%
% 6.1 %
\section{Almost complex structures in $\R^4$}
\markboth{6.1. Almost complex structures in $\R^4$}
{6.1. Almost complex structures in $\R^4$}

\hspace{13.5pt}
Let us consider an almost complex structure $j$ in $\R^4$ in a neighborhood
of zero. Suppose that $(N_j)_0\ne0$. In this case the image
$\op{Im}N_j(\cdot,\cdot)$ is two-dimensional and complex generated by a
vector $N_j(\x,\e)\ne0$ for any two complex independent vectors $\x$ and $\e$.
Consider a two-dimensional distribution
$\Pi^2_x=\op{Im}(N_j)_x\subset T_x=T_x\R^4$. This distribution is an
invariant of almost complex structure and hence invariants of this
distribution lead to invariants of the almost complex structure.

An important
invariant of a distribution is the Tanaka invariant (\cite{T}, \cite{Y}):
a collection of graded Lie algebras ${\scriptstyle\cal Q}(x)$ for each point
$x$, associated to the filtered algebras $\{D^p(x)\subset T_x\}_{p\ge1}$,
where the module of sections ${\cal D}^{(p)}$ of the distribution $D^p$ is
defined as $(p-1)$-th derivative: ${\cal D}^{(p)}=\p^{(p-1)}{{\cal D}^{(1)}}=
\p^{(p-2)}{{\cal D}^{(1)}}+ [{\cal D}^{(1)}, \p^{(p-2)}{\cal D}^{(1)}]$,
$[\cdot,\cdot]$ being the commutator of the vector fields, and
$\p^{(0)}{\cal D}^{(1)}={\cal D}^{(1)}$ is the sections module of the initial
distribution $D^1$; the Lie product is induced by the commutator of
vector fields.

Let us suppose that the first derivative of the distribution $\Pi^2_*$ is
nontrivial at the point: $(\p^{(1)}\Pi^2)_0\ne \Pi_0^2$. Then we have a
three-dimensional distribution $\Pi^3 = \p^{(1)}\Pi^2$ in $\ok(0)$.
There exist vectors $\x_i=\x_i(x)\ne0$, $i=1, 2, 3$, in the space
$\Pi^3_x\subset T_x$ such that $N_j(\x_1,\x_3)=\x_1$, $N_j(\x_2,\x_3)=-\x_2$.
Actually, for any $\x_3\in\Pi^3\setminus\Pi^2$ the mapping $\e\mapsto
N(\e,\x_3)$ is an orientation reversing isomorphism of $\Pi^2$, and hence
there exist two one-dimensional invariant subspaces. Consider a vector $\x_1$
on one of them. We have: $N(\x_1,\x_3)=f\x_1$, $f\ne0$. Let's change the
vector $\x_3\mapsto \dfrac1f\x_3$. Then $N(\x_1,\x_3)=\x_1$ and
$N(j\x_1,\x_3)=-j\x_1$, i.e. we may take $\x_2=j\x_1$.

Denote by $\U_1= <\x_1>$ and $\U_2= <\x_2>$ the one-dimensional subspaces
generated by the vectors $\x_1$ and $\x_2$, $\Pi^2_x=(\U_1)_x\oplus(\U_2)_x$.
Let us note that if the half-space $(\Pi_x^3)^+\subset \Pi^3_x\setminus
\Pi^2_x$, in which lies the vector $\x_3$, is fixed then the vector is
defined up to $\Pi^2_x$-shifts: $\x_3\mapsto\x_3+\a_1\x_1+\a_2\x_2$. If
one changes the half-space: $\x_3\mapsto{\tilde\x}_3=-\x_3$, then the
vectors $\x_1$ and $\x_2$ (and hence the distributions $\U_1$ and $\U_2$)
interchanges:
${\tilde\x}_1=\pm\x_2$, ${\tilde\x}_2=\pm\x_1$. Selecting a half-space
in $\Pi^3_x\setminus\Pi^2_x$, i.e. an orientation on $\T^1_x=\Pi^3_x/\Pi_x^2$,
we fix which one of the distributions $\U_1$ and $\U_2$ is the first and
which is the second; changing this orientation we change the numeration.
In other words, there is a canonical orientation on the two-dimensional space
$\T^1_x\t P(\Pi^2_x)$. Let us call this orientation the $\T$-orientation.
Note also that there's singled out a pair of two-dimensional affine subspaces
in the space $\Pi^3_x$: $\{\pm\x_3+\a_1\x_1+\a_2\x_2\}$, i.e. there's fixed
a metric on $\T^1_x$. Let us call it the $\T$-metric.

Consider now the cofactor $\Xi^1_x= T_x^4/\Pi_x^3$ of the subspace
$\Pi^3_x\subset T_x$. There exists a vector $\x_4\notin\Pi_x^3$ such that
$N_j(\x_1,\x_4)=\x_2$, $N_j(\x_2,\x_4)=-\x_1$. In its half-space
$T^4_x\setminus\Pi^3_x$ it is determined up to $\Pi^2_x$-shifts:
$\x_4\mapsto\x_4+\a_1\x_1+\a_2\x_2$, and so there exists a natural metric on
$\Xi^1_x$. Let us call it the $\Xi$-metric.
Under the change of half-space $\x_4\mapsto{\tilde\x}_4=-\x_4$
in one of two subspaces $\U_1$ or $\U_2$ the orientation changes:
$\x_1\mapsto -\x_1$ или $\x_2\mapsto-\x_2$. Therefore there's determined
the natural orientation of the space $\U_1^1\t\U_2^1\t\Xi^1$. Let us
call it the $\Xi$-orientation.

Let us call the $\U\T\Xi$-invariant the pair of distribution $\U_i$ and
$\T$- and $\Xi$-metrics and orientations. From the arguments above
together with theorem~6 it follows

 \begin{th}\po
Let the Nijenhuis tensors of almost complex structures $j_1$ in
$\ok_{\R^4}(x)$ and $j_2$ in $\ok_{\R^4}(y)$ do not equal to zero and the
first derivative of the distribution $\Pi^2$ is nontrivial,
$\partial^{(1)}\Pi^2\ne\Pi^2$. Pseudoholomorphic or antipseudoholomorphic
mapping $u: (\ok_{\R^4}(x),j_1) \to (\ok_{\R^4}(y),j_2)$ exists if and only
if there exists a mapping $\ok_{\R^4}(x) \to \ok_{\R^4}(y)$, which
transforms one $\U\T\Xi$-invariant to the other.
\qed
 \end{th}

 \begin{rk}{6}
In the case if the second derivative of the distribution $\Pi^2$ does not
coincide with the first, $\partial^{(2)}\Pi^2\ne\partial^{(1)}\Pi^2\ne\Pi^2$,
i.e. $(\partial^{(2)}\Pi^2)_x=T_x^4$, the Tanaka invariant (\cite{T})
--- graded Lie algebra --- has the underlying space of the form
 $$
{\scriptstyle\cal Q}(x)=
{\scriptstyle\cal Q}(x)_1\oplus
{\scriptstyle\cal Q}(x)_2\oplus
{\scriptstyle\cal Q}(x)_3 \simeq \Pi^2_x\oplus
\T^1_x\oplus \Xi^1_x.
 $$
In this
case, because of the gradation, the Lie product is determined by a 2-form on
$\Pi^2_x$ with values in $\T^1_x$ and by 1-form on $\Pi^2_x$ with values in
$\op{Hom}(\T^1_x,\Xi^1_x)$. Due to existence of the $\T$- and $\Xi$-metrics
and orientations and due to the decomposition
$\Pi^2_x=(\U_1)_x\oplus(\U_2)_x$, the Lie algebra structure on
${\scriptstyle\cal Q}(x)$
is given by elements $\oo^2_x\in (\U_1)^*_x\otimes(\U_2)^*_x=
\op{Hom}(\U_1,\U_2^*)_x$ and $\oo^1_x\in (\U_1)^*_x\oplus(\U_2)^*_x$. Thus
the Tanaka invariant gives us two additional invariants: $\oo^1$ and
$\oo^2$. \qed
 \end{rk}

%%%%%%%%%%%%%%%%%%%%%%%%%%%%%%%%%%%%%%%%%%%%%%%%%%%%%%%%%%%%%%%%%%%%%%%%%%
% 6.2 %
\section{Lie almost complex structures}
\markboth{6.2. Lie almost complex structures}
{6.2. Lie almost complex structures}

\hspace{13.5pt}
Let us consider an almost complex structure $j$ and its Nijenhuis tensor
$N_j$. This tensor defines a bilinear skew-symmetric map
$N_j: {\cal D}\otimes {\cal D}\to {\cal D}$ on the module ${\cal D}$ of vector
fields. In general this map is not a Lie product. Actually, from the Jacobi
identity for the elements $\x$, $\e$, $j\zeta$
 $$
N(\x,N(\e,j\zeta))+ N(\e,N(j\zeta,\x))+ N(j\zeta,N(\x,\e))=0
 $$
it follows that
 $$
jN(\x,N(\e,\zeta))+jN(\e,N(\zeta,\x))- jN(\zeta,N(\x,\e))=0,
 $$
which together with the Jacobi identity for $\x$, $\e$, $\zeta$ gives:
 \begin{equation}
N(\x,N(\e,\zeta))\equiv0.
% \tag{25}
\add
 \end{equation}

Thus $({\cal D}, N(\cdot,\cdot))$ is a Lie algebra iff identity (25) holds.
Let us denote by $\Pi_x\subset T_x$ the subspace generated by
$\op{Im} N(\cdot,\cdot)$ and denote $A_x=\{\x\,\vert\, N(\x,\cdot)\equiv0\}$.
Formula (25) is equivalent to the inclusion of ($j$-invariant) subbundles
$\Pi\subset A\subset TX$.

 \begin{dfn}{6}
Let us call an almost complex structure $j$ on $X$ Lie almost complex
structure if its Nijenhuis tensor $N_j$ defines a Lie algebra structure
on the module~${\cal D}$.
 \end{dfn}

From theorem 6 (chapter~5) and the arguments above it follows

 \begin{th}\po
A Lie almost complex structure $j$ defines a solvable Lie algebra ${\cal G}
=({\cal D},N_j(\cdot,\cdot))$ with the solvability rank~2: ${\cal G}^{(2)}=0$.
Lie almost complex structures of general position are equivalent or
antiequivalent (formally or smooth) iff the corresponding Lie algebras
$({\cal D}, N(\cdot,\cdot))$ are equivalent. \hfill $\Box$
 \end{th}

Thus the theorem reduces the classification of Lie almost complex structures
to the classification of Lie algebras bundle (Lie algebras structures on the
tangent spaces $T_{x^\prime}X$), i.e. to the analysis of Lie multiplication
deformations; let us note that not all such deformations realize, see
remark~3.

 \begin{rk}{7}
Let us consider {\it Nijenhuis (sub)algebra\/} $\A^\N_j$, which is
generated by the complex structure $j$ by means of the algebraic commutator
$[A,B]=B\barwedge A- (-1)^{\op{deg}A\op{deg}B} A\barwedge B$ and
the Nijenhuis bracket $\lv\cdot,\cdot\rv$, in the algebra of vector-valued
forms $\O^\bullet\otimes{\cal D}$ on the manifold $X$ (see~\cite{FN},
\cite{ALV}, the notations are similar to~\cite{KS}). Here $\O^\bullet$ is
the algebra of scalar valued differential forms and ${\cal D}$ is the
module of vector fields. Note that in the case of Lie almost complex
structure $\A^\N_j=$ $<j, N_j, j\c N_j>$. Actually, as formulae
(2.9), (2.10), (2.16), (5.6), (5.8), (5.22), (5.23) from \cite{FN} together
with the Jacobi identity $N_j\barwedge N_j=0$ for $N_j$ show, the following
equalities hold true:
 \begin{eqnarray*}
&&[j,j]=0,\ [j,N_j]=-3j\barwedge N_j=\frac32 N_j\barwedge j=-3j\c N_j,\
[j,j\c N_j]=3N_j,\\
&&[N_j,N_j]=2N_j\barwedge N_j=0,\ [N_j,j\c N_j]=0,\ [j\c N_j,j\c N_j]=0,\\
&&\lv j,j\rv=N_j,\ \lv j, N_j\rv=0,\ \lv j,j\c N_j\rv=-\frac12 N_j\barwedge
N_j=0,\\
&&\lv N_j,N_j\rv=0,\ \lv N_j,j\c N_j\rv=0,\ \lv j\c N_j,j\c N_j\rv=0.
 \end{eqnarray*}
It follows from the example below that the values set of $(\A^\N_j)_x$ at the
point $x$ is not sufficient for formal classification of the almost complex
structure, i.e. $\A^\N_j\ne\A^\infty_j$, see the beginning of~chapter~2. \qed
 \end{rk}

 \begin{ex}{6}
Let us consider the almost complex structure
defined by the formulae:
 \begin{eqnarray*}
j\p_1=\p_2,{\ \ \!}&\qquad j\p_3=\p_4+f(x^5)\p_1,{\ \,}&\qquad j\p_5=\p_6,\\
j\p_2=-\p_1,     & \qquad j\p_4=-\p_3-f(x^5)\p_2,      & \qquad j\p_6=-\p_5.
 \end{eqnarray*}
The Nijenhuis tensor on the basis vectors is given by the identities:
 $$
N_j(\p_1,\cdot)=N_j(\p_2,\cdot)=N_j(\p_3,\p_4)=N_j(\p_5,\p_6)=0,
 $$
 $$
N_j(\p_3,\p_5)=-N_j(\p_4,\p_6)=f^\prime(x^5)\p_2,\
N_j(\p_3,\p_6)=N_j(\p_4,\p_5)=f^\prime(x^5)\p_1.
 $$
One can see that the structure $j$ is not formally classified by the algebra
$(\A^\N_j)_x$. Actually, for $f(x^5)=\ve\cdot x^5$, there is an additional
independent invariant $(dN_j)_0$, and for $f(x^5)=x^5+\ve(x^5)^2$
the algebra $\A^\N_j$ also is not sufficient for distinguishing the
structures corresponding to different $\ve$, even though the Nijenhuis
tensors are conjugated in all the point from the neighborhood of $x$. \qed
 \end{ex}

 \begin{rk}{8}
Let us consider the Nijenhuis tensor of a Lie almost complex structure
$j$ in $\ok(x)$. There exists an increasing $j$-invariant filtration
$\T_k\subset T_x$, where $\T_k$ is generated by
$\op{Im}d^kN(\cdot,\cdot,\ast)$. One can relate to a point $x$ the
following invariant: the associated graded ($j$-)complex Lie algebra
regardness $N$-product --- $\oplus V_k$, $V_k=\T_k/\T_{k-1}$. \qed
 \end{rk}

\appendix
%%%%%%%%%%%%%%%%%%%%%%%%%%%%%%%%%%%%%%%%%%%%%%%%%%%%%%%%%%%%%%%%%%%%%%%%%%
% app %
 \chapter{Linear Nijenhuis tensors of general position}
\markboth{General position}{General position}

\hspace{13.5pt}
Here we prove that tensors $N\in\N_j$ satisfying the condition of
definition~4 of chapter~4 are actually the tensors of general
position in $\N_j$. Let $V\simeq\C^n$ be the complex linear space.

 \begin{th}\po
The set $\u_j$ of all tensors $N\in\N_j(V)$ such that $\n_N(\x)=2$ for almost
all $\x\in V$ is open and everywhere dense in $\N_j$.
 \end{th}

 \begin{proof}
$1^\circ$. Let us fix a complex basis on $V$, i.e. we fix an isomorphism
$V=\C^n$. Let us also fix a Euclidean structure on $V$ in which the basis
is orthonormal; in particular the automorphism $j$ is an orthogonal
transformation. Consider an $(n-1)$-dimensional complex subspace
$\Pi\subset V$, i.e. $(2n-2)$-dimensional subspace invariant under $j$.
Consider the mapping $N_\x:\Pi\to V$, $N_\x(\e)=N(\x,\e)$, where
$\x\notin\Pi$. The last condition means that the complex line $\C\x=
<\x,j\x>$ is transversal to $\Pi$. Let us fix an orthonormal complex basis
$\e_1,\dots,\e_{2n-2}$ in $\Pi$, $\e_{2k}=j\e_{2k-1}$. Define the function
$\a_N: V\setminus\Pi\to\R$ as the $(2n-2)$-dimensional volume of the body
$\{\sum\limits_{k=1}^{2n-2} x_k\cdot N_\x(\e_k)\ \vert\ 0\le x_k\le1\}$.
The condition $\n_N(\x)=2$ of definition~4 is equivalent to $\a_N(\x)>0$
for $\x\in V\setminus\Pi$.

Note that the square of the function $\a_N(\x)$ is a rational fracture of the
variable $\x=(\x_1,\dots,\x_{2n})$. So if $\a_N\not\equiv0$
then the set of $\x$ such that $\a_N(\x)=0$ is a hypersurface (actually a
cone over a surface) and so has the measure~0.

$2^\circ$. Show that the set $\u_j$ in $\N_j$ is nonempty. Really, consider
the tensor $N\in\N_j$ defined by
 $$
N(\vec z,\vec w)=\left(0,
\left\vert
\begin{array}{cc}
\bar{z_1} & \bar{z_2} \\
\bar{w_1} & \bar{w_2}
\end{array}
\right\vert,
\left\vert
\begin{array}{cc}
\bar{z_1} & \bar{z_3} \\
\bar{w_1} & \bar{w_3}
\end{array}
\right\vert,
\dots,
\left\vert
\begin{array}{cc}
\bar{z_1} & \bar{z_n} \\
\bar{w_1} & \bar{w_n}
\end{array}
\right\vert
\right),
 $$
where $\vec z=(z_1,\dots,z_n)\in\C^n$, $\vec w=(w_1,\dots,w_n)\in\C^n$ and
$\bar z_k=x_k-jy_k$ for $z_k=x_k+jy_k$. One easily sees that if $z_1\ne0$
and $\vec w\ne\lambda\vec z$ for any $\lambda\in\C$ then $N(\vec z,\vec
w)\ne0$. So here we have $\a_N(\x)\ne0$ for all $\x\notin\Pi=\{\vec z \
\vert\ z_1=0\}$.

$3^\circ$. Show that $\u_j$ is open in $\N_j$. Let $N\in\u_j\subset\N_j$.
Consider the set $\{\x\ \vert\ \a_N(\x)=0,\ \vert\x\vert=1\}$. This is
a hypersurface in $S^{2n-1}\subset V$. Let $\w$ be a small neighborhood of
it. Let $a>0$ be the minimum of $\a_N$ on the closure
$\overline{S^{2n-1}\setminus\w}$. There exists a neighborhood $\Y$ of $N$
in $\N_j$ such that for any $N'\in\Y$ we have
 $$
\a_{N'}(\x)>\frac a2,\ \x\in\overline{S^{2n-1}\setminus\w}.
 $$
Hence $\a_{N'}\not\equiv0$ and $N'\in\u_j$ according to $1^\circ$, and we
have proved the desired inclusion $\Y\subset\u_j$.

$4^\circ$. Let us prove that $\u_j$ is everywhere dense in $\N_j$. Let
$N\in\N_j$ and let $N_0\in\u_j$ be the tensor constructed in $2^\circ$.
Consider the 1-parameter family $N_\ve=\ve N+ (1-\ve)N_0 \in\N_j$. We have:
$N_1=N$.  The square of the function $\a_{N_\ve}(\x)$ is a rational
fracture of $\ve$. This function is not
equal to zero for $\ve=0$ for almost all $\x$. So it could not be zero
identically. In particular there exists a number $\ve_0>0$ such that for
any $\ve\in (1-\ve_0,1)$ the function $\a_{N_\ve}$ does not vanish
identically, i.e. $N_\ve\in\u_j$. \qed
 \end{proof}

In the same way one could prove the theorem "without $j$". Let the set $\V$
consist of all vector two-forms $A\in\wedge^2(\R^n)^*\otimes\R^n$ such that
for almost all $\x\in\R^n$ the equality $A(\x,\e)=0$ implies that
$\e=\lambda\x$ for some $\lambda\in\R$.

 \begin{th}\po
$\V$ is open and everywhere dense in $\wedge^2(\R^n)^*\otimes\R^n$. \qed
 \end{th}

One of the forms $A\in\V$ is given in the example~3 from chapter~4. We
prolong this tensor $A$ to the tensor $N$ on the complexified space
$\R^n_\C=\C^n\simeq\R^n\oplus\R^n=\R^n\oplus j\R^n$ as is shown in the
example. Let us show that $\n_N=2$ almost everywhere. Really, consider
$\x_1=(1,\dots,1)$, $\e_1=(1,\dots,1)$,
$\Pi=\{(0,*,\dots,*)\oplus(0,*,\dots,*)\}\subset\R^n\oplus j\R^n.$
Let $\x_2=(0,\x_2^2,\dots,\x_2^n)$ and $\e_2=(0,\e_2^2,\dots,\e_2^n)$ be
two vectors from $\Pi$ ($\x_2^1=\e_2^1=0$). We have:

 \begin{eqnarray*}
N(\x_1\oplus\e_1,\x_2\oplus\e_2)=&&\!\!\!\!\!\!\!
\left[\sum_{k=1}^n(\x_2^{k+1}-\x_2^k)e_k-\sum_{k=1}^n(\e_2^{k+1}-\e_2^k)e_k
\right]\\
\oplus&&\!\!\!\!\!\!\!
\left[-\sum_{k=1}^n(\e_2^{k+1}-\e_2^k)e_k-\sum_{k=1}^n(\x_2^{k+1}-\x_2^k)e_k
\right].
 \end{eqnarray*}

So if $N_{\x_1\oplus\e_1}(\x_2\oplus\e_2)=0$ then
$\sum\limits_{k=1}^n(\x_2^{k+1}-\x_2^k)e_k=0$,
$\sum\limits_{k=1}^n(\e_2^{k+1}-\e_2^k)e_k=0$,
from where we have $\x_2^k=\e_2^k=0$, $k=2,\dots,n$.
Hence $\a_N(\x_1\oplus\e_1)>0$ and $\n_N(\x_1\oplus\e_1)=2$, which means
that $N$ is of general position (see the proof of theorem~9).

 \begin{rk}{9}
We can construct tensors $A\in\V$, $A:\wedge^2\R^n\to\R^n$, such that
$A(\x,\e)=0$ only for parallel vectors $\x$ and $\e$. Let us denote them by
$\V_0$. For $n=2$: $\V_0=[\wedge^2(\R^n)^*\otimes\R^n]\setminus\{0\}$. For
$n=3$ we can take as example of $A$ the usual vector product
$\R^3\wedge\R^3\to\R^3$. Let us consider the case $n=4$. We have the
Pl\"ucker-Grassmann mapping:
 $$
PG: \wedge^2\R^4\to\R^4,\ (\x.\e)\mapsto \te=(\te_{rs})_{1\le r<s\le4},
 $$
where $\te_{rs}=\left\vert\begin{array}{cc} \x^r&\x^s \\
                                            \e^r&\e^s \end{array}\right\vert$.
 The image $\op{Im}(PG)$ is the quadric (cone) given by the equation
 $$
\CC:\ \{\te_{12}\te_{34}-\te_{13}\te_{24}+\te_{14}\te_{23}=0\}.
 $$

Let us take the four-dimensional subspace $\{\te_{12}=-\te_{34},\
\te_{13}=\te_{24}\}$ in $\R^6$. The orthogonal
projection $\op{pr}$ on this subspace is given by the formulae:
 \begin{eqnarray*}
\mu_1 &\!\!\!=\!\!\!& \te_{12}-\te_{34}, \\
\mu_2 &\!\!\!=\!\!\!& \te_{13}+\te_{24}, \\
\mu_3 &\!\!\!=\!\!\!& \te_{14}, \\
\mu_4 &\!\!\!=\!\!\!& \te_{23}.
 \end{eqnarray*}
So by composition $\op{pr}\c PG$ we obtain the map $A:\wedge^2\R^4\to\R^4$.
This map satisfies the property that $A(\x,\e)=0$ only for parallel $\x$
and $\e$. Actually, in this case $\te_{12}=\te_{34}$, $\te_{13}=-\te_{24}$,
$\te_{14}=0$, $\te_{23}=0$. Substituting these equalities into the formula
for $\CC$ we have
 $$
\te_{12}^2+\te_{13}^2=0,
 $$
i.e. $\te_{rs}=0$ for all $1\le r<s\le4$, which means that $\x$ is parallel
to $\e$.

As in the formula of example~3 we could construct a mapping $N_A:
\wedge^2\R^{2n}\to\R^{2n}$, $\R^{2n}=\R^n\oplus\R^n=\R^n\oplus j\R^n$,
$N_A\in\N_j$. But this map does not satisfy the property that
$N_A(\x,\e)=0$ only for $\C$-dependent vectors $\x$ and $\e$ (note that in
difference with definition~4 here we require that the property holds not
almost always but always). Actually, the author do not know whether there
is for any $n$ the tensor $N\in\N_j^{(n)}$ such that
$N(\x,\e)=0$ implies $\x\in\C\e$ unless $\e=0$. \qed
 \end{rk}
%\newpage
%
%
%
%%%%%%%%%%%%%%%%%%%%%%%%%%%%%%%%%%%%%%%%%%%%%%%%%%%%%%%%%%%%%%%%%%%%%%%%%%
%
%%%%%%%%%%%%%%%%%%%%%%%%%%%%%%%%%%%%%%%%%%%%%%%%%%%%%%%%%%%%%%%%%%%%%%%%%%
%

\bigskip
\bigskip
\bigskip
\bigskip
\bigskip
\bigskip
\bigskip
\bigskip
\bigskip

{\it \hspace{-19pt} Address:}
{\footnotesize
 \begin{itemize}
  \item
P. Box 546, 119618, Moscow, Russia
  \item
Chair of Applied Mathematics, Moscow State Technological University
\linebreak
{\rm n. a.} Baumann, Moscow, Russia
 \end{itemize}
}

{\it \hspace{-19pt} E-mail:} \quad
{\footnotesize
lychagin\verb"@"glas.apc.org or borkru\verb"@"difgeo.math.msu.su
}

%%%%%%%%%%%%%%%%%%%%%%%%%%%%%%%%%%%%%%%%%%%%%%%%%%%%%%%%%%%%%%%%%%%%%%%%%%
\end{document}